\documentclass[11pt,reqno]{amsart}

\usepackage{amsmath,amssymb,amsfonts,amsthm,amscd,mathtools,latexsym,enumerate}
\usepackage[
  colorlinks=true,
  linkcolor=red,
  citecolor=red,
  urlcolor=red
]{hyperref}
\usepackage{xcolor}
\usepackage{mathrsfs}

\newcommand{\aut}{\mathrm{Aut}}
\newcommand{\id}{\mathrm{id}}
\newcommand{\zz}{\mathbb{Z}}
\newcommand{\Ftwo}{\mathbb{F}_2}
\newcommand{\Dcov}{\mathrm{D}}
\newcommand{\gp}{\mathrm{GP}}
\newcommand{\autp}{\aut^{+}}
\newcommand{\autm}{\aut^{-}}
\newcommand{\ol}{\overline}
\newcommand{\bip}{\mathbin{\star}}
\newcommand{\Bipprod}{\mathop{\star}}
\newcommand{\Und}{\mathrm{U}}
\newcommand{\op}{\mathrm{op}}
\newcommand{\mb}{\mathbf{Bip}_2}
\newcommand{\ole}{\overline{\boldsymbol{\epsilon}}}
\newcommand{\mfw}{\mathfrak{W}}
\newcommand{\ostab}{\operatorname{Stab}}
\newcommand{\rse}{\operatorname{Res}}
\newcommand{\mrc}{ \varrho}

\newtheorem{theorem}{Theorem}[section]
\newtheorem{lemma}[theorem]{Lemma}
\newtheorem{corollary}[theorem]{Corollary}
\newtheorem{proposition}[theorem]{Proposition}

\theoremstyle{definition}
\newtheorem{definition}[theorem]{Definition}
\newtheorem{remark}[theorem]{Remark}
\newtheorem{example}[theorem]{Example}
\newtheorem{problem}[theorem]{Problem}

\begin{document}

\title[Fibre factorization and automorphisms of bipartite products]{Fibre factorization and automorphisms of direct products of bipartite graphs}

\author{Xiaomeng Wang}
\address{School of Mathematics and Statistics, Lanzhou University, Lanzhou, Gansu 730000, P.\,R. China}
\email{wangxiaomeng@lzu.edu.cn}


\date{}

\begin{abstract}
We determine the automorphism groups of direct products of finite connected twin-free bipartite graphs whose factors are different from $K_2$. The product is generally disconnected, so the problem requires both the
automorphism groups of its components and a criterion for component isomorphism. We retain the order of the bipartition classes and use the fibre product over $K_2$. An edge-relation construction, together with McKenzie's strict refinement theorem, yields unique fibre-prime
factorization and a coordinate description of every specified isomorphism. After refining the original factors, component isomorphism classes are
described by intersections of factor-symmetry orbits with a set of admissible labels. The stabilizer of a label acts faithfully as the fullautomorphism group of its component. These data determine the full automorphism group, including products with repeated or fibre-composite factors. For products of fibre powers of one prime graph, we compute the component groups and express their multiplicities through a subset-sum
polynomial. We also obtain an explicit semidirect product for a mixed product with one nonbipartite factor and explain its relation to known stability results. An application recovers the known double-cover groups
of bipartite generalized Petersen graphs.
\end{abstract}

\renewcommand{\subjclassname}{2020 Mathematics Subject Classification}

\maketitle

\section{Introduction}\label{sec:introduction}
\subsection{Motivation}
The automorphism problem for a graph product asks how its symmetries
are determined by the factors. For a direct product of connected bipartite
graphs, there are two questions: what is the automorphism group of each
connected component, and which components are isomorphic? The first
determines the symmetries within components; the second determines which
components may be permuted. We answer both questions from prime-factor
data for connected twin-free bipartite factors, none isomorphic to $K_2$.
The factorization used here is for the fibre product over $K_2$.

All graphs considered in this paper are finite and undirected and have no multiple edges.
We write $\Gamma_0$ for the class of graphs in which loops are allowed, otherwise graphs are loopless.
Write $x\sim_X y$ when $x$ and $y$ are adjacent in $X$. For a graph $X$, let $V(X)$, $E(X)$, and $\aut(X)$ denote its vertex set, edge set, and automorphism group, respectively,  and let $N_X(x)$ denote the neighbourhood of a vertex $x$. A graph is \emph{twin-free} if distinct vertices have distinct neighbourhoods.
We follow the graph terminology of~\cite{Big01,HIK2011}, the group terminology of~\cite{Hun00} and the categorical terminology of~\cite{Mac1998}.

The \emph{direct product} \cite{HIK2011} of two graphs $X$ and $Y$ in $\Gamma_0$, denoted by $X\times Y$, has vertex set $V(X)\times V(Y)$, where $(x,y)$ is adjacent to $(x',y')$ if and only if $\{x,x'\}\in E(X)$ and $\{y,y'\}\in E(Y)$.
This product is commutative and associative up to the natural coordinate isomorphisms, and its unit is the graph with one vertex and one loop.
A nonunit graph is \emph{prime} if at least one factor is the unit in every factorization as a direct product. Two graphs are \emph{coprime} if they have no common direct factor of order greater than one.
For nonempty graphs $X$ and $Y$, the coordinatewise action gives a natural embedding
\begin{equation}\label{eq:direct}
\aut(X)\times\aut(Y)\leq\aut(X\times Y).
\end{equation}
Determining when equality holds in~\eqref{eq:direct}, and describing the additional automorphisms when the inclusion is strict, are basic problems in the theory of graph products.
A complete answer must account for the way automorphisms interact with the factor coordinates.

For connected nonbipartite graphs, prime factorization supplies the relevant coordinates.
Such graphs have unique prime factorizations with respect to the direct product in $\Gamma_0$~\cite{HIK2011, I1998, McK71}. Under the additional twin-free hypothesis, every isomorphism between products of prime factors is induced by factor isomorphisms and a permutation of the corresponding coordinates~\cite{D1974,HIK2011}.
Thus prime factorization determines not only the isomorphism types of the factors but also the coordinate form of isomorphisms between products.

For bipartite factors, two different issues intervene. First, Weichsel's theorem~\cite{Weichsel1962}, also stated in~\cite[Theorem~5.9]{HIK2011}, implies that a direct product of $r\geq1$ connected bipartite graphs, each with at least one edge, has  $2^{r-1}$ connected components. A description of its full automorphism group therefore requires both the automorphism groups of individual components and a criterion for deciding which components are isomorphic. Second, prime direct product factorizations of connected bipartite graphs need not be unique~\cite{HIK2011,HP2015}. Thus the nonbipartite coordinate theory cannot simply be applied to the components of the direct product of bipartite graphs.
The component isomorphism question has a precise answer for two factors. Jha, Klav\v{z}ar and Zmazek~\cite{JKZ1997} proved that the two components are isomorphic if one factor admits an automorphism interchanging its bipartition classes, and conjectured the converse. Hammack studied the problem in~\cite{H2006} and proved the converse in~\cite{H2009}. A separate result of Hammack and Puffenberger~\cite{HP2015}, extending~\cite{H2013}, shows that the unique bipartite prime factor  occurring in any prime factorization of a connected bipartite graph is determined up to isomorphism. That theorem concerns a connected direct product  with exactly one bipartite prime factor. Here we study the connected components of a product in which  every factor is bipartite, and we seek the structure of their isomorphisms as well as their isomorphism types.

Our approach retains an ordering of the two bipartition classes and replaces the direct product by the fibre product over $K_2$.  For connected factors, the fibre product is the component in which the coordinate colours agree. It remains connected and retains the factor coordinates.  We encode such ordered bipartite graphs by binary relations on their edge sets, prove that the relevant class of relations is closed under direct factors, and use McKenzie's refinement theory~\cite{McK71} to obtain fibre prime factorizations and coordinate descriptions of specified isomorphisms. The resulting graph-theoretic description keeps track of factorization, bipartition reversals, and permutations of repeated isomorphic factors simultaneously.

The remaining problem is to recover the component isomorphism classes and their automorphism groups from prime-factor data. The wreath-product formula for the automorphism group of a disjoint union is classical. The essential task is to identify its constituent groups and their multiplicities directly from the factors, without treating the components as independent input data. This leads to the following question.

\begin{problem}\label{pro:pro1}
 How  can one describe the automorphism groups of the connected components of direct products of connected bipartite graphs and account for the interaction between bipartition reversals and permutations of repeated isomorphic factors?
    \end{problem}

\subsection{Main results}
The results proceed in three stages. First, Theorem~\ref{th:thuf} identifies the fibre-prime factors. Next, Theorems~\ref{th:thci} and~\ref{th:thfc} show how a specified isomorphism acts on their coordinates, including a common colour-reversal sign. Finally, Theorem~\ref{th:thrfa} uses these data to determine component isomorphism classes, their automorphism groups, and the automorphism group of the full direct product. 

We first introduce the notation needed for the statements. Write $\Ftwo$ for the field with two elements, $\mathbf0_t,\mathbf1_t\in\Ftwo^t$ for the constant vectors, and $Q_t=\Ftwo^t/\langle\mathbf1_t\rangle$ for $t\geq1$.
A bar denotes the class of a binary vector in this quotient. Write $S_t$ for the symmetric group on $\{1,\ldots,t\}$, with $S_0=1$. Maps act on the left and are composed from right to left. We use the same convention for relations: $x(S\circ R)z$ means that $xRy$ and $ySz$ for some $y$. Write $\zz_n=\mathbb Z/n\mathbb Z$, and denote the path and cycle on $n$ vertices by $P_n$ and $C_n$, respectively. For a finite set $J$, write $S_J$ for its symmetric group. For a group $G$, the notation $G\wr S_t=G^t\rtimes S_t$ uses the action $(\sigma\cdot\mathbf{g})_i=g_{\sigma^{-1}(i)}$; we set $G\wr S_0=1$. Bold lower-case letters denote tuples, whereas their subscripted entries are individual coordinates.

An \emph{ordered bipartite graph} is a pair $(B,\mrc_B)$, where $B$ is a bipartite graph and $\mrc_B:B\to K_2$ is a surjective graph homomorphism. Equivalently, $B$ has a fixed ordered bipartition $V(B)=B^0\sqcup B^1$, where $B^i=\mrc_B^{-1}(i)$ for $i\in\{0,1\}$. We call $\mrc_B$ the \emph{colour map} and write $\Und(B)$ for the underlying graph of $B$, obtained by forgetting the order of the two classes. A morphism from $B$ to $C$ is a graph homomorphism $f:\Und(B)\to\Und(C)$ satisfying
\begin{equation*}
 \mrc_C\circ f=\mrc_B.
\end{equation*}

Ordered bipartite graphs and their colour-preserving homomorphisms form a category, denoted by $\mb$.
A graph isomorphism $f:\Und(B)\to \Und(C)$ is called an \emph{ordinary isomorphism}. Such an isomorphism need not preserve the colour classes.

\begin{definition}\label{def:biproduct}
Let $B$ and $C$ be ordered bipartite graphs. Their \emph{fibre product over $K_2$}, denoted by $B\bip C$, is the ordered bipartite graph with vertex set
\[
 V(B\bip C)
 =\{(b,c)\in V(B)\times V(C):\mrc_B(b)=\mrc_C(c)\}
\]
and edge set
\[
 E(B\bip C)=\{\{(b,c),(b',c')\}:\;(b,c),(b',c')\in V(B\bip C),
 \{b,b'\}\in E(B), \{c,c'\}\in E(C)\}.
\]
The colour map is $\mrc_{B\bip C}(b,c)=\mrc_B(b)=\mrc_C(c)$, so its two ordered colour classes are
\[
 (B\bip C)^0=B^0\times C^0,
 \qquad
 (B\bip C)^1=B^1\times C^1.
\]
The graph $K_2$, equipped with its identity colour map, is the unit of this product. A connected nonunit ordered bipartite graph $B$ is called \emph{fibre-prime} (or, \emph{$f$-prime}) if, in every factorization $B\cong X\bip Y$ in $\mb$, at least one of $X$ and $Y$ is isomorphic to $K_2$.
\end{definition}


The fibre product is associative and commutative up to the natural coordinate isomorphisms. We write $\Bipprod_{i=1}^r B_i$ for an iterated fibre product and take the empty product to be $K_2$. The fibre product has a different unit and notion of primality from the ordinary direct product. The first theorem gives a common refinement of any two fibre product decompositions and, in particular, unique factorization into nonunit $f$-prime factors.

\begin{theorem}\label{th:thuf}
Let $B$ be a finite connected twin-free ordered bipartite graph. Let $r,s\geq0$, with empty fibre products interpreted as $K_2$, and suppose
\[
 B\cong G_1\bip\cdots\bip G_r
 \cong H_1\bip\cdots\bip H_s
\]
in $\mb$.
Then the following statements hold.
\begin{enumerate}[(i)]
\item\label{it:ituf1} There are connected twin-free $f$-prime ordered bipartite graphs $B_\lambda$ and two indexed partitions of a finite set $\mathcal S$ (with empty parts allowed),
\[
 \mathcal S=U_1\sqcup\cdots\sqcup U_r
 =V_1\sqcup\cdots\sqcup V_s
\]
such that
\[
 B\cong\Bipprod_{\lambda\in \mathcal S}B_\lambda,
 \qquad
 G_i\cong\Bipprod_{\lambda\in U_i}B_\lambda,
 \qquad
 H_j\cong\Bipprod_{\lambda\in V_j}B_\lambda.
\]
Thus $B$ admits a factorization into $f$-prime factors.

\item\label{it:ituf2} If $G_1\bip\cdots\bip G_r$ and $H_1\bip\cdots\bip H_s$ are already factorizations into nonunit $f$-prime factors, then $r=s$ and there is a permutation $\sigma\in S_r$ such that for each $1\leq i\leq r$, $G_i\cong H_{\sigma(i)}$ in $\mb$.
Consequently, the multiset of nonunit $f$-prime factors of $B$ is unique up to isomorphism in $\mb$.
\end{enumerate}
\end{theorem}

Unique factorization identifies the prime factors up to colour-preserving isomorphism and permutation. Determining automorphism groups requires additional information, namely a coordinate description of each specified isomorphism. The next theorem provides this description, covering cases where isomorphic factors repeat.

\begin{theorem}\label{th:thci}
Let $B=B_1\bip\cdots\bip B_r$ and $C=C_1\bip\cdots\bip C_s$
be products of connected twin-free nonunit $f$-prime ordered bipartite graphs, where $r,s\geq1$. If $\varphi:B\to C$ is an isomorphism in $\mb$, then $r=s$, and there exist a permutation $\sigma\in S_r$ and isomorphisms $\varphi_i:B_i\longrightarrow C_{\sigma(i)}$ in $\mb$ for $1\leq i\leq r$
such that, whenever $\varphi(b_1,\ldots,b_r)=(c_1,\ldots,c_r)$, we have
\begin{equation}\label{eq:eqcis}
 c_{\sigma(i)}=\varphi_i(b_i),
\end{equation}
for every $1\leq i\leq r$.
\end{theorem}

An ordinary isomorphism between connected bipartite graphs either preserves both bipartition classes or interchanges them. For a fibre product, these choices cannot be made independently in different coordinates since all factor maps must have the same colour sign. The next theorem gives this precise compatibility condition. It is the  fibre product counterpart of the coordinate description for connected nonbipartite twin-free  direct products in~\cite[Theorem~8.15]{HIK2011}.

\begin{theorem}\label{th:thfc}
Let $B$ and $C$ be connected twin-free ordered bipartite graphs with fixed nonunit $f$-prime decompositions $B\cong B_1\bip\cdots\bip B_r$ and $C\cong C_1\bip\cdots\bip C_s$ in $\mb$, where $r,s\geq1$.
Let $\varphi:\Und(B)\longrightarrow\Und(C)$ be an ordinary graph isomorphism. Then $r=s$. Moreover, there are a unique global sign $\delta\in\Ftwo$, a permutation $\sigma\in S_r$, and ordinary isomorphisms
$\varphi_i:\Und(B_i)\longrightarrow\Und(C_{\sigma(i)})$,  for $1\leq i\leq r$, such that, for each $1\leq i\leq r$ and every $x\in V(B_i)$,
\begin{equation}\label{eq:eqfas}
 \mrc_{C_{\sigma(i)}}(\varphi_i(x))
 =\mrc_{B_i}(x)+\delta.
\end{equation}
 If $\varphi(b_1,\ldots,b_r)=(c_1,\ldots,c_r)$,
then for each $1\leq i\leq r$,
\begin{equation}\label{eq:eqocis}
 c_{\sigma(i)}=\varphi_i(b_i).
\end{equation}
Thus all corresponding factor isomorphisms preserve their ordered bipartitions when $\delta=0$, and all reverse them when $\delta=1$.
\end{theorem}

We now return to the  direct product. Its components correspond to choices of factor colourings, with two choices identified when every colour is reversed. Theorem~\ref{th:thfc} will turn isomorphisms between these components into permutations and isomorphisms of fibre prime factors. We first describe the component labels and the action of factor symmetries on them.
For an ordered bipartite graph $B$, let $B^{\op}$ have the same underlying graph with the two colour classes interchanged, and put $B^{(0)}=B$ and $B^{(1)}=B^{\op}$. Given connected ordered bipartite graphs $B_1,\ldots,B_r$ and $\boldsymbol{\epsilon}=(\epsilon_1,\ldots,\epsilon_r)\in\Ftwo^r$, the fibre product $U=B_1^{(\epsilon_1)}\bip\cdots\bip B_r^{(\epsilon_r)}$ has a colour map given, for ${\bf x}=(x_1,\ldots,x_r)\in V(U)$, by
\begin{equation}\label{eq:eqgcd} \mrc_U({\bf x})=\mrc_{B_1}(x_1)+\epsilon_1=\cdots=\mrc_{B_r}(x_r)+\epsilon_r\in\Ftwo.
\end{equation}
We call $\mrc_U({\bf x})$ the \emph{global colour} of ${\bf x}$. The two sets
$U^0=\mrc_U^{-1}(0)$ and $U^1=\mrc_U^{-1}(1)$ are the \emph{global bipartition classes}.

Reversing all factor colourings leaves the underlying graph of $U$ unchanged and reverses only its global colouring. Thus component labels lie in the quotient $Q_r$ defined above.
For $\boldsymbol{ \epsilon}=(\epsilon_1,\ldots,\epsilon_r)\in\Ftwo^r$, write $\ol{\boldsymbol{\epsilon}}$ for its class in $Q_r$ and define
\[ \Omega_{\ol{\boldsymbol{\epsilon}}}
 =\Und\bigl(B_1^{(\epsilon_1)}\bip\cdots\bip B_r^{(\epsilon_r)}\bigr),\qquad P=\Und(B_1)\times\cdots\times\Und(B_r).\]
For ${\bf x}=(x_1,x_2,\ldots,x_r)\in V(P)$,  its \emph{coordinate colour vector} is $\pi({\bf x})=(\mrc_{B_1}(x_1),\ldots,\mrc_{B_r}(x_r))$. Proposition~\ref{pro:procd} shows that
${\bf x}\in V(\Omega_{\ol{\boldsymbol{\epsilon}}})$ precisely when $\overline{\pi({\bf x})}=\ole$, and that the graphs $\Omega_{\ol{\boldsymbol{\epsilon}}}$ are exactly the $2^{r-1}$ connected components of $P$. 

Choose $\sigma\in S_r$ and ordinary isomorphisms $f_i:\Und(B_i)\to\Und(B_{\sigma(i)})$ for every $i$. With $\mathbf{f}=(f_1,\ldots,f_r)$, define $w_{\sigma,\mathbf{f}}:P\to P$ by $w_{\sigma,\mathbf{f}}({\bf x})={\bf y}$ for ${\bf x}=(x_i)_{i=1}^r$ and ${\bf y}=(y_i)_{i=1}^r$, where, for $1\leq i\leq r$,
\begin{equation*}
 y_{\sigma(i)}=f_i(x_i).
\end{equation*}
Equivalently, $[w_{\sigma,\mathbf{f}}({\bf x})]_{\sigma(i)}
 =f_i(x_i)$,
where $[{\bf z}]_j$ denotes the $j$-th coordinate of a tuple ${\bf z}$. Thus $f_i$ acts on the  $i$-th input coordinate, and $\sigma$ specifies its output position. Each $w_{\sigma,\mathbf{f}}$ is bijective and preserves adjacency, and its inverse is determined by $\sigma^{-1}$ and the maps $f_i^{-1}$. These maps form the \emph{factor symmetry group} $\mfw=\mfw(B_1,\ldots,B_r)$.
For clarity, the composition rule is $w_{\tau,\mathbf g}w_{\sigma,\mathbf f}=w_{\tau\sigma,\mathbf h}$,
where $h_i=g_{\sigma(i)}\circ f_i$.
Since each $B_i$, $1\leq i\leq r$, is connected, each $f_i$ has a unique sign $s_i\in\Ftwo$, determined by the identity
 \begin{equation*}
 \mrc_{B_{\sigma(i)}}(f_i(x_i))=\mrc_{B_i}(x_i)+s_i.
\end{equation*}
Here $x_i$ ranges over $V(B_i)$. Put ${\bf s}=(s_1,\ldots,s_r)$ and define $\sigma_*{\bf v}$ by $(\sigma_*{\bf v})_{\sigma(i)}=v_i$ for ${\bf v}\in\Ftwo^r$. The induced action on component indices is
\begin{equation}\label{eq:eqfsa}
 w_{\sigma,\mathbf{f}}\cdot\ol{\boldsymbol{\epsilon}}
 =\overline{\sigma_*(\boldsymbol{\epsilon}+{\bf s})}.
\end{equation}
This is well defined since $\sigma_*{\bf 1}_r={\bf 1}_r$, where ${\bf 1}_r=(1,\ldots,1)$. Write $\ostab_{\mfw}(\ol{\boldsymbol{\epsilon}})=\{w\in\mfw:w\cdot\ol{\boldsymbol{\epsilon}}=\ol{\boldsymbol{\epsilon}}\}$ for the stabilizer of the indicated quotient class.

The original factors $B_i$ need not be $f$-prime, so their coordinate symmetries alone need not describe all component isomorphisms. We therefore refine each $B_i$ into its $f$-prime factors. The original coordinates become blocks of refined coordinates, and the components belonging to the original product are selected by requiring the refined colour labels to be constant within each block.

Let $r\geq1$, and let $B_1,\ldots,B_r$ be connected twin-free nonunit ordered bipartite graphs. By Theorem~\ref{th:thuf}, for each $1\leq i\leq r$, fix an isomorphism in $\mb$
\begin{equation}\label{eq:eqcr}
 \theta_i:B_i\longrightarrow\Bipprod_{k\in I_i}C_k,
\end{equation}
where $I_1,\ldots,I_r$ partition $\{1,\ldots,N\}$ into nonempty sets, and every $C_k$ is connected, twin-free, nonunit, and $f$-prime. The order of the coordinates within each $I_i$ is fixed. Write $\mathcal{I}=(I_1,\ldots,I_r)$, and define
\[ \widetilde P=\prod_{k=1}^N\Und(C_k),\qquad
 \widetilde{\mfw}=\mfw(C_1,\ldots,C_N).\]
We use the same quotient notation $Q_N$ for the refined coordinates.
For a vector $\mathbf{a}=(a_1,\ldots,a_N)\in\Ftwo^N$, define
\[
 \Lambda_{\ol{\mathbf{a}}}=\Und(C_1^{(a_1)}\bip\cdots\bip C_N^{(a_N)}).
\]
Each $C_k$ is a factor graph, whereas $\Lambda_{\ol{\mathbf {a}}}$ is a connected component of $\widetilde{P}$. An element $\mathbf{x}=(x_1,\ldots,x_N)\in V(\Lambda_{\ol{\mathbf{a}}})$ is a vertex, with $x_k\in V(C_k)$ for each $k$.

Define the repetition map $d:\Ftwo^r\to\Ftwo^N$ and its induced map $\ol d:Q_r\to Q_N$ by
\[
 (d(\boldsymbol{\epsilon}))_k=\epsilon_i,\quad \ol d(\ol{\boldsymbol{\epsilon}})=\ol{d(\boldsymbol{\epsilon})},
\]
where $\boldsymbol{\epsilon}= (\epsilon_1,\ldots,\epsilon_r)\in \Ftwo^r$ and $k\in I_i$.
Since $d(\mathbf{1}_r)=\mathbf{1}_N$, the map $\ol{d}$ is well defined. Write
\begin{equation*}
 \mathcal L_{\mathcal I}=\ol{d}(Q_r)=\{\ol{\mathbf a}\in Q_N:a_k=a_\ell\text{ whenever }k,\ell\in I_i \text{ for some }i\}.
\end{equation*}
The condition on the right is unchanged when $\mathbf a$ is replaced by $\mathbf a+\mathbf{1}_N$.
For $\mathbf x=(x_1,\ldots,x_r)\in V(P)$, write $\theta_i(x_i)=(u_k)_{k\in I_i}$ for each $i$, and define $\Theta(\mathbf x)=(u_1,\ldots,u_N)\in V(\widetilde P)$.
We write $\Theta_{\ol{\boldsymbol{\epsilon}}}$ for  the restriction $\Theta|_{V(\Omega_{\ol{\boldsymbol{\epsilon}}})}$ (its image will be identified in Theorem~\ref{th:thrfa}(\ref{it:ierfaa})).

The group $\widetilde{\mfw}$ acts on $Q_N$ by Eq.~\eqref{eq:eqfsa}, now using the $C_k$-coordinates. The subset $\mathcal L_{\mathcal I}$ need not be invariant; all orbits below are taken in the ambient label space $Q_N$. Let $\mathscr{O}_{\mathcal{I}}$ be the set of its orbits that meet $\mathcal L_{\mathcal I}$. For each $\mathcal{O}\in\mathscr{O}_{\mathcal{I}}$,  choose $\ol{\mathbf a}_{\mathcal{O}}\in\mathcal{O}\cap \mathcal L_{\mathcal I}$ and let $n_{\mathcal{O}}=|\mathcal{O}\cap \mathcal L_{\mathcal I}|$.

\begin{theorem}\label{th:thrfa}
Let $r\geq1$, and let $B_1,\ldots,B_r$ be connected twin-free nonunit ordered bipartite graphs. Put
\[
 P=\Und(B_1)\times\cdots\times\Und(B_r).
\]
For each $1\leq i\leq r$, fix an $f$-prime factorization
$\theta_i:B_i\longrightarrow\Bipprod_{k\in I_i}C_k$
in $\mb$ as in~\eqref{eq:eqcr}.
With the notation introduced above for these fixed factorizations, the following statements hold.
\begin{enumerate}[(i)]
\item\label{it:ierfaa} The map $\Theta$ is an isomorphism $P\longrightarrow\sqcup_{\ol{\mathbf{a}}\in \mathcal L_{\mathcal I}}\Lambda_{\ol{\mathbf{a}}}$,
and $\Theta_{\ol{\boldsymbol{\epsilon}}}:\Omega_{\ol{\boldsymbol{\epsilon}}}\to\Lambda_{\ol{d}(\ol{\boldsymbol{\epsilon}})}$ is an isomorphism for every $\ol{\boldsymbol{\epsilon}}\in Q_r$. Moreover, $|\mathcal L_{\mathcal I}|=2^{r-1}$.
\item\label{it:ierfab} For $\ol{\boldsymbol{\epsilon}},\ol{\boldsymbol{\eta}}\in Q_r$, the components $\Omega_{\ol{\boldsymbol{\epsilon}}}$ and $\Omega_{\ol{\boldsymbol{\eta}}}$ are isomorphic if and only if $\ol{d}(\ol{\boldsymbol{\epsilon}})$ and $\ol{d}(\ol{\boldsymbol{\eta}})$ belong to the same $\widetilde{\mfw}$-orbit in $Q_N$.
\item\label{it:ierfac}  For every $\ol{\boldsymbol{\epsilon}}\in Q_r$, restriction followed by conjugation by $\Theta_{\ol{\boldsymbol{\epsilon}}}$ defines an isomorphism from $\ostab_{\widetilde{\mfw}}(\ol d(\ol{\boldsymbol{\epsilon}}))$ to $\aut(\Omega_{\ol{\boldsymbol{\epsilon}}})$ mapping elements $w\in\ostab_{\widetilde{\mfw}}(\ol d(\ol{\boldsymbol{\epsilon}}))$ to $\Theta_{\ol{\boldsymbol{\epsilon}}}^{-1}\circ
     (w|_{\Lambda_{\ol d(\ol{\boldsymbol{\epsilon}})}})\circ\Theta_{\ol{\boldsymbol{\epsilon}}}$.
\item\label{it:ierfad}
The automorphism group of $P$ is given by
\begin{equation}\label{eq:eqgau}
 \aut(P)\cong\prod_{\mathcal{O}\in\mathscr{O}_{\mathcal{I}}}
 \left(\ostab_{\widetilde{\mfw}}(\ol{\mathbf a}_{\mathcal{O}})
                  \wr S_{n_{\mathcal{O}}}\right).
\end{equation}
\item\label{it:ierfae}
Let $T_1,\ldots,T_q$ represent the distinct isomorphism types of the underlying graphs among the $C_k$. For each $1\leq \lambda\leq q$, let $J_\lambda=\{k:\Und(C_k)\cong T_\lambda\}$ and $t_\lambda=|J_\lambda|$. Then
\begin{equation*}
 \widetilde{\mfw}\cong\prod_{\lambda=1}^q
       \left(\aut(T_\lambda)\wr S_{t_\lambda}\right).
\end{equation*}
\end{enumerate}
\end{theorem}

Theorem~\ref{th:thrfa} answers Problem~\ref{pro:pro1} for connected twin-free nonunit factors. Parts~(\ref{it:ierfab}) and~(\ref{it:ierfac}) show that  two selected components are isomorphic precisely when their refined labels lie in the same orbit, and the full automorphism group of a selected component is the stabilizer of its label. Part~(\ref{it:ierfad}) then assembles the independent component automorphisms and the permutations of isomorphic components. Part~(\ref{it:ierfae}) expresses the refined factor symmetry group using the automorphism groups and multiplicities of the underlying $f$-prime factors.

\subsection{Significance and relation to previous work}
The essential content of Theorem~\ref{th:thrfa} is the determination of component data from the refined factors. The wreath product formula for a disjoint union is standard. Here its constituent groups and multiplicities are obtained as stabilizers and orbit-intersection sizes.
The component groups and component isomorphism types are not additional inputs. In particular, specified isomorphisms between components, rather than only the existence of such isomorphisms, must be controlled.

The dependence on classical theory is explicit. McKenzie's strict refinement theorem~\cite{McK71} is the decomposition input. Lemmas~\ref{le:legrc} and~\ref{le:lefc} supply the edge-relation correspondence and its closure under factors. These lemmas transfer specified relation isomorphisms back to graph coordinates. Theorem~\ref{th:thfc} then handles the common colour sign, and the original blocks specify the components selected in Theorem~\ref{th:thrfa}. The correspondence is used for objects and isomorphisms (no assertion about all noninvertible homomorphisms is needed).
This theory concerns fibre products, not uniqueness of ordinary direct
prime factorizations of bipartite graphs.

Hammack's two-factor theorem~\cite{H2009} determines when the two components are isomorphic. Our treatment applies to any number of bipartite factors and recovers both the component isomorphism classes and the component automorphism groups via a common action on refined factor labels. The restriction to connected twin-free factors is essential throughout. We describe all component isomorphisms and automorphism groups for arbitrary numbers of factors. Corollary~\ref{co:corp} makes the description explicit for products of fibre powers of one prime graph.  This calculation uses the refinement of the original factors and distinguishes the present setting from the singleton-block specialization.

Following~\cite{QXZZ2021}, a graph pair $(X,Y)$ is \emph{stable} if equality holds in \eqref{eq:direct}. The \emph{canonical double cover} is $\Dcov(X)=X\times K_2$, and $X$ is \emph{stable} when $(X,K_2)$ is stable.
See~\cite{FH2022,LMS2015,MSZ1989,QXZ2019,Surowski2001,Surowski2003,Wilson2008} for graph stability and \cite{GLX2025,QXZ2024,QXZZ2021,WQX2025,WXZ2025,WG2026} for graph pairs. Corollary~\ref{pr:prmp} computes the mixed-product group when $\Dcov(X)$ and the bipartite factor are nonisomorphic fibre-prime graphs.
Its stability equivalence is covered by the factor-loopless theorem in~\cite[Theorem~1.3]{WG2026}. We explain the implication between the
hypotheses after the proof. The final generalized Petersen formulas recover~\cite[Theorem~1.2(iii)]{QXZ2021} and are not new classifications.


\subsection{Organization of the paper}
Section~\ref{sec:secpre} introduces the necessary notation and basic results on ordered bipartite graphs, fibre products, and binary relation structures. In Section~\ref{sec:biproduct}, we prove Theorems \ref{th:thuf}, \ref{th:thci} and \ref{th:thfc}, and derive the component isomorphism criterion in Corollary~\ref{cor:corsfc}.
In Section~\ref{subsec:automorphisms}, we prove Theorem~\ref{th:thrfa} and derive Corollaries~\ref{co:corfa} and~\ref{th:thadp} as special
cases. We then compute products of fibre powers in Corollary~\ref{co:corp}, derive the semidirect product formula of Corollary~\ref{pr:prmp} and obtain the automorphism group of canonical double covers of bipartite generalized Petersen graphs in Corollary \ref{co:coqxzbb}.

\section{Preliminaries}\label{sec:secpre}
We introduce notation for ordered bipartite graphs, establish basic properties of their fibre products, and describe their correspondence with binary relation structures.
\subsection{Ordered bipartite graphs}

Recall that $B^{(\epsilon)}$ denotes a recolouring of $B$, whereas $B^\epsilon$ is one of its vertex classes. In this subsection $B$ and $C$ are connected ordered bipartite graphs. An isomorphism in $\mb$ preserves the two ordered classes; an ordinary isomorphism may interchange them.
For any ordinary isomorphism $f:\Und(B)\rightarrow\Und(C)$, the uniqueness of the bipartition up to interchange implies that there is a unique sign $\chi(f)\in\Ftwo$ such that, for every vertex $x\in V(B)$,
\begin{equation}\label{eq:eqis}
 \mrc_C(f(x))=\mrc_B(x)+\chi(f).
\end{equation}
The map $f$ preserves the two colour classes when $\chi(f)=0$ and interchanges them when $\chi(f)=1$.

For a connected ordered bipartite graph $B$, define $\chi_B:\aut(\Und(B))\to\Ftwo$ by $\chi_B(g)=\chi(g)$ for $g\in\aut(\Und(B))$. Applying Eq.~\eqref{eq:eqis} twice, we obtain $\chi_B(gh)=\chi_B(g)+\chi_B(h)$ for $g,h\in\aut(\Und(B))$. Thus $\chi_B$ is a group homomorphism.
Define
\begin{align*}
 \autp(B)&=\ker\chi_B
 =\{g\in\aut(\Und(B)):g(B^0)=B^0\},\\
 \autm(B)&=\chi_B^{-1}(1)
 =\{g\in\aut(\Und(B)):g(B^0)=B^1\}.
\end{align*}
The automorphism group of $B$ as an object of $\mb$ is precisely $\autp(B)$. If an automorphism interchanges the two colour classes, then $\autp(B)$ has index two in $\aut(\Und(B))$, and $\autm(B)$ is its nontrivial coset. In all cases, $\aut(\Und(B))=\autp(B)\sqcup\autm(B)$. If no automorphism interchanges the colour classes, then $\autm(B)=\varnothing$.

For connected ordered bipartite graphs $B$ and $C$, define the group \[\mathcal{G}(B,C)=\{(g,h)\in\aut(\Und(B))\times\aut(\Und(C)):\chi_B(g)=\chi_C(h)\}.\]
In other words, \[\mathcal{G}(B,C)=\bigl(\autp(B)\times\autp(C)\bigr)\cup\bigl(\autm(B)\times\autm(C)\bigr),\]
where the second set is empty if either factor has no automorphism interchanging its two colour classes. Since $\mathcal{G}(B,C)$ is the kernel of the homomorphism $(g,h)\mapsto\chi_{{B}}(g)+\chi_{C}(h)$ from $\aut(\Und(B))\times\aut(\Und(C))$ to $\Ftwo$, it is a subgroup of $\aut(\Und(B))\times\aut(\Und(C))$.

 For the canonical double cover $\Dcov(X)$ defined in the introduction,  we use the ordered bipartition
\[
 \Dcov(X)^0=V(X)\times\{0\},
 \qquad
 \Dcov(X)^1=V(X)\times\{1\},
\]
and write $\kappa_X((x,i))=(x,1-i)$ for the covering involution, where $(x,i)\in V(\Dcov(X))$. Recall that a graph $X$ is called \emph{stable} if $\aut(\Dcov(X))=\aut(X)\times\langle\kappa_X\rangle$ \cite{MSZ1989}.

\begin{lemma}\label{le:leosdo}
Let $B$ be a connected twin-free ordered bipartite graph. If all vertices in one colour class have degree one, then $B\cong K_2$ in $\mb$. Consequently, if $B$ is nonunit, each colour class has at least two vertices and contains a vertex of degree at least two.
\end{lemma}

\begin{proof}
Suppose that every vertex of $B^0$ has degree one. If $B^1$ contained two vertices, a path between them would have an internal vertex in $B^0$ of degree at least two. Thus $B^1$ has one vertex. Since $B$ is connected and twin-free $B\cong K_2$. The same argument applies with the colours reversed.

If either colour class has only one vertex, then $B$ is a star by the connectedness, and moreover, $B\cong K_2$ since $B$ is twin-free. Consequently, a nonunit graph has at least two vertices in each class. Since $B$ is connected, each class contains a vertex of degree at least two.
\end{proof}

\subsection{The fibre product over $K_2$ }
We first establish some basic properties of the fibre product over $K_2$. Throughout this subsection, we assume that every ordered bipartite graph has two nonempty colour classes and at least one edge.

Let $B, C$ and $D$ be ordered bipartite graphs. The natural maps
\[
 (B\bip C)\bip D\longrightarrow B\bip(C\bip D),
 \qquad
 ((b,c),d)\longmapsto(b,(c,d)),
\]
and
\[
 B\bip C\longrightarrow C\bip B,
 \qquad
 (b,c)\longmapsto(c,b),
\]
are isomorphisms in $\mb$, where $b\in V(B)$, $c\in V(C)$ and $d\in V(D)$. Hence the fibre product is associative and commutative up to these natural isomorphisms. We write $\Bipprod_{i=1}^r B_i$.
For connected factors, $B\bip C$ is connected, whereas $\Und(B)\times \Und(C)$ has two connected components. Thus the fibre product is precisely the component in which the coordinate colours agree. Its unit and its notion of primality differ from  those of the direct product.

We now establish the fibre product properties used in the sequel.

\begin{lemma}\label{le:lebp}
Let $B$ and $C$  be connected ordered bipartite graphs. Then $B\bip C$ is connected.  
\end{lemma}

\begin{proof}
 Let $(b,c)$ and $(b',c')$ be vertices of $B\bip C$. Then $\mrc_B(b)=\mrc_C(c)$ and $\mrc_B(b')=\mrc_C(c')$. Choose a path $P:b=b_0,b_1,\ldots,b_p=b'$
in $B$ and a path $Q:c=c_0,c_1,\ldots,c_q=c'$
in $C$. We have $p\equiv \mrc_B(b')-\mrc_B(b)\pmod{2}$ and $q\equiv \mrc_C(c')-\mrc_C(c)\pmod{2}$, so $p\equiv q\pmod{2}$.
Without loss of generality, assume that $p\leq q$. Then $q-p$ is even. Since $B$ is connected and has an edge, every vertex of $B$ has a neighbour. We obtain a walk $b=b_0',b_1',\ldots,b_q'=b'$ of length $q$ by inserting the walk $z,z',z$  at a vertex $z$ of $P$ exactly $(q-p)/2$ times, where $z'$ is a neighbour of $z$. At every step, both coordinates move along an edge, and their colours change simultaneously. Therefore $(b_0',c_0),(b_1',c_1),\ldots,(b_q',c_q)$
is a walk in $B\bip C$ from $(b,c)$ to $(b',c')$. Every walk contains a path with the same endpoints, so $B\bip C$ is connected.
\end{proof}

The next result indexes all connected components of a direct product of bipartite graphs by  coordinate colour vectors.

\begin{proposition}\label{pro:procd}
Let $B_1,\ldots,B_r$ be connected ordered bipartite graphs, where $r\geq2$ is an integer, and let $Q_r=\Ftwo^r/\langle{\bf 1_r}\rangle$. Let $P=\Und(B_1)\times\cdots\times\Und(B_r)$.
Then the following statements hold.
\begin{enumerate}[(i)]
\item\label{it:itocd1} The graphs $\Omega_{\ole}$, indexed by $\ole\in Q_r$, are exactly the connected components of  $P$.
\item\label{it:itocd2} $P$ has $2^{r-1}$ connected components.
\item\label{it:itocd4} The representatives $\boldsymbol{\epsilon}$ and $\boldsymbol{\epsilon}+{\bf 1}_r$ determine the same underlying component with opposite ordered bipartitions.
\end{enumerate}
\end{proposition}

\begin{proof}
(\ref{it:itocd1}). For any vertex ${\bf x}=(x_1,\ldots,x_r)$ in $V(P)$, its coordinate colour vector is $\pi({\bf x})=\bigl(\mrc_{B_1}(x_1),\ldots,\mrc_{B_r}(x_r)\bigr)\in\Ftwo^r$.
If ${\bf x}$ is adjacent to ${\bf y}=(y_1,\ldots,y_r)$ in $P$, then $x_i$ is adjacent to $y_i$ in $B_i$ for every $1\leq i\leq r$. Since every edge of a bipartite graph joins vertices of opposite colours, we have $\mrc_{B_i}(y_i)=\mrc_{B_i}(x_i)+1$ for every $1\leq i\leq r$. Thus $\pi({\bf y})=\pi({\bf x})+{\bf 1}_r$.
The class $\overline{\pi({\bf x})}\in Q_r$ is constant on every edge of $P$ and is therefore constant on every connected component of $P$.

Fix $\boldsymbol\epsilon=(\epsilon_1,\ldots,\epsilon_r)\in\Ftwo^r$. By the definition of the fibre product, a vertex ${\bf x}=(x_1,\ldots,x_r)$ belongs to $
 V\bigl(B_1^{(\epsilon_1)}\bip\cdots\bip B_r^{(\epsilon_r)}\bigr)$
if and only if $\mrc_{B_i}(x_i)+\epsilon_i=\mrc_{B_j}(x_j)+\epsilon_j$ for every $1\leq i, j\leq r$. Thus there is a value $t\in\Ftwo$ such that  $\mrc_{B_i}(x_i)=\epsilon_i+t$ for every $1\leq i\leq r$.
By Eq.~\eqref{eq:eqgcd},  $t$ is the global colour of ${\bf x}$ in the fibre product. Hence $\pi({\bf x})=\boldsymbol\epsilon+t{\bf 1}_r$,
which is equivalent to $\overline{\pi({\bf x})}=\overline{\boldsymbol{\epsilon}}$ in $Q_r$.
Therefore, the vertex set of $\Omega_{\ol{\boldsymbol{\epsilon}}}$ is precisely the set of vertices of $P$ whose coordinate colour vectors have class $\ol{\boldsymbol{\epsilon}}$.

By Lemma~\ref{le:lebp}, each fibre product $B_1^{(\epsilon_1)}\bip\cdots\bip B_r^{(\epsilon_r)}$
is connected.
Since no edge joins a vertex in this set to a vertex whose colour vector lies in a different class, $\Omega_{\ol{\boldsymbol{\epsilon}}}$ is a connected component of $P$.

(\ref{it:itocd2}).
Each factor has two nonempty colour classes.
Thus every vector in $\Ftwo^r$ occurs as a coordinate colour vector and  every quotient class  indexes a component.
As $Q_r=\Ftwo^r/\langle\mathbf1_r\rangle$ has order
$2^{r-1}$, the number of components is  $2^{r-1}$.

(\ref{it:itocd4}).
 The vectors  $\boldsymbol\epsilon$ and
$\boldsymbol\epsilon+\mathbf1_r$  represent the same class in $Q_r$ and therefore specify  the same underlying component. Replacing one representative by the other changes the common colour $t$ in the preceding display to  $t+1$, which interchanges the two global colour classes of that component.
\end{proof}

\begin{lemma}\label{le:lebpa}
Let $B$ and $C$  be connected ordered bipartite graphs. Then $B\bip C$ is twin-free if and only if both $B$ and $C$ are twin-free. Moreover, $\Und(B)\times\Und(C)
 =\Und(B\bip C)\sqcup\Und(B\bip C^{\op})$.
\end{lemma}

\begin{proof}
By Proposition \ref{pro:procd}(\ref{it:itocd1}), $B\bip C$ is a connected component of $\Und(B)\times\Und(C)$. Together with \cite[Lemma~2.3]{HIK2011}, this shows that $B\bip C$ is twin-free if and only if both $B$ and $C$ are twin-free.
Moreover, when $r=2$, the two elements of $Q_2$ are represented by $(0,0)$ and $(0,1)$.
By Proposition~\ref{pro:procd}, we have $\Und(B)\times\Und(C)
 =\Und(B\bip C)\sqcup\Und(B\bip C^{\op})$.
\end{proof}

We next give a useful sufficient condition for fibre primality.

\begin{lemma}\label{le:lecp}
Let $B$ be a connected twin-free nonunit ordered bipartite graph. Then $B$ is $f$-prime if one of the following statements holds:
\begin{enumerate}[(i)]
\item\label{it:itcpa} all vertices of one colour class have the same prime degree;
\item\label{it:itcpb} $\Und(B)$ contains no cycle of length four.
\end{enumerate}
\end{lemma}

\begin{proof}
Suppose that $B\cong X\bip Y$ in $\mb$, where both factors are nonunit. Since both colour classes of each factor are nonempty, the coordinate projections are surjective graph homomorphisms. Hence $X$ and $Y$ are
connected. By Lemma~\ref{le:lebpa}, since $B$ is twin-free, $X$ and $Y$ are twin-free. For $a\in X^i$, $c\in Y^i$ and $i\in\{0,1\}$,
\[
 \deg_{X\bip Y}(a,c)=\deg_X(a)\deg_Y(c).
\]
Under (\ref{it:itcpa}), interchange the colour labels if necessary so that the common prime degree occurs on $B^0$, and denote it by $p$. Fix $a_0\in X^0,c_0\in Y^0$. Varying $a$ in $\deg_X(a)\deg_Y(c_0)=p$ makes $\deg_X$ constant on $X^0$. Varying $c$ with $a_0$ fixed makes $\deg_Y$ constant on $Y^0$. The two positive constants multiply to the prime $p$, so $\deg_X$ is $1$ or $\deg_Y$  is $1$.
Lemma~\ref{le:leosdo} makes that factor $K_2$, a contradiction.

Under (\ref{it:itcpb}), Lemma~\ref{le:leosdo} supplies $a\in X^0$ with distinct neighbours $a_1,a_2\in X^1$, and $c\in Y^1$ with distinct neighbours $c_1,c_2\in Y^0$. The four vertices
\[
 (a,c_1),\ (a_1,c),\ (a,c_2),\ (a_2,c)
\]
are distinct, and form a four-cycle in $X\bip Y$, again a contradiction.
\end{proof}

\subsection{Relation structures on ordered bipartite graphs}\label{subsec:relations}

We now construct a binary relation structure on the edge set of a connected twin-free ordered bipartite graph. We begin with definitions from the theory of relational structures~\cite{Ch1967,McK71}.

\begin{definition}\cite{McK71}\label{de:derds}
A \emph{relation structure} is a pair $\mathcal{X}=(X,R)$, where $X$ is a nonempty set and $R\subseteq X\times X$ is a binary relation. For $x\in X$, write $R[x]=\{y\in X:xRy\}$ and $R^{-1}[x]=\{y\in X:yRx\}$. The structure is \emph{reflexive} if $xRx$ for every $x\in X$. It is \emph{connected} if for any two elements $x, y\in X$, there are an integer $m\geq 0$ and elements $x=x_0, x_1,\ldots,x_m=y$ such that for every $0\leq i\leq m-1$, either $x_iRx_{i+1}$ or $x_{i}R^{-1}x_{i+1}$.

Define two equivalence relations on $X$ by
\begin{equation}\label{eq:eqcp}
 x\equiv_0x'
 \quad\Longleftrightarrow\quad
 R[x]=R[x'],
 \qquad
 x\equiv_1x'
 \quad\Longleftrightarrow\quad
 R^{-1}[x]=R^{-1}[x'].
\end{equation}
The relation structure is \emph{thin} if
$(\equiv_0\cap\equiv_1)=\Delta_X$, where $\Delta_X=\{(x,x):x\in X\}$ is the diagonal relation. We call it \emph{rectangular} in this paper if
\begin{equation}\label{eq:eqrectang}
 xRy
 \quad\Longleftrightarrow\quad
 \text{there exists }z\in X\text{ such that }z\equiv_0x\text{ and }z\equiv_1y.
\end{equation}
\end{definition}
For any $x\in X$, write $[x]_{\equiv_0}=\{y\in X:x\equiv_0y\}$ and $[x]_{\equiv_1}=\{y\in X:x\equiv_1y\}$ for its two equivalence classes.
Write $X/\equiv_0=\{[x]_{\equiv_0}: x\in X\}$ and
$X/\equiv_1=\{[x]_{\equiv_1}: x\in X\}$.
For relation structures $\mathcal{X}=(X,R)$ and $\mathcal{Y}=(Y,Q)$,
a relation isomorphism is a bijection $f:X\to Y$ satisfying
$xRx'\iff f(x)Qf(x')$. Their direct product is $(X\times Y, R\times Q)$, where $(x,y)(R\times Q)(x',y')$ if and only if $xRx'$ and $yQy'$.
When two underlying sets occur in one formula, superscripts $R,Q$, or $R\times Q$
on $\equiv_i$ identify the relation defining the equivalence class.
The one-element reflexive relation structure is the unit. A nonunit relation structure is \emph{$d$-prime} if, in every factorization into two relation structures, at least one factor is isomorphic to the unit.
Denote by $\mathfrak{G}$ the class of all finite connected reflexive thin rectangular relation structures, together with the one-element unit.

Let $B$ be a connected twin-free ordered bipartite graph with $V(B)=B^0\sqcup B^1$.
We identify an edge with the ordered pair $(a,b)$ of its endpoints, where $a\in B^0$ and $b\in B^1$. For edges $e=(a,b)$ and $f=(c,d)$, define
\[
 e\,\rho_0\,f\quad\Longleftrightarrow\quad a=c,
 \qquad
 e\,\rho_1\,f\quad\Longleftrightarrow\quad b=d.
\]
Let $R_B=\rho_0\circ\rho_1$ and $\mathcal{R}(B)=(E(B),R_B)$.
With the right-to-left convention fixed above, $eR_Bf$ means that there exists an edge $g$ of $B$ such that $e\rho_0g$ and $g\rho_1f$. Thus
\begin{equation}\label{eq:eqcer}
 (a,b)R_B(c,d)
 \quad\Longleftrightarrow\quad
 (a,d)\in E(B).
\end{equation}

Conversely, let $\mathcal{X}=(X,R)\in\mathfrak{G}$. Define an ordered bipartite graph $\mathsf{B}(\mathcal{X})$ with $\mathsf{B}(\mathcal{X})^0=(X/\equiv_0)\times\{0\}$ and $\mathsf{B}(\mathcal{X})^1=(X/\equiv_1)\times\{1\}$.
For every $x\in X$, include the edge $e_x=\{([x]_{\equiv_0},0),([x]_{\equiv_1},1)\}$.
If $e_x=e_y$, then $x\equiv_0y$ and $x\equiv_1y$. Since $\mathcal{X}$ is thin, we obtain $x=y$. Thus, distinct elements of $X$  correspond to distinct edges, so $\mathsf{B}(\mathcal{X})$ is a simple bipartite graph.

The next lemma shows that the two constructions recover one another and respect the corresponding products.

\begin{lemma}\label{le:legrc}
Let $B$ be a connected twin-free ordered bipartite graph, and let $\mathcal{X}=(X,R),\mathcal{Y}=(Y,Q)\in\mathfrak{G}$. Then the following statements hold.
\begin{enumerate}[(i)]
\item\label{it:itgrc1} The relation structure $\mathcal{R}(B)=(E(B),R_B)$ belongs to $\mathfrak{G}$, and the two equivalence relations defined in  \eqref{eq:eqcp} are exactly $\rho_0$ and $\rho_1$.
\item\label{it:itgrc2}
The bipartite graph $\mathsf{B}(\mathcal{X})$ is connected and twin-free, and $\mathcal{R}(\mathsf{B}(\mathcal{X}))\cong\mathcal{X}$.
Every relation isomorphism $\mathcal{X}\to\mathcal{Y}$ induces an isomorphism in $\mb$ from $\mathsf{B}(\mathcal{X})$ to $\mathsf{B}(\mathcal{Y})$.
\item\label{it:itgrc3}
For connected twin-free ordered bipartite graphs $B$ and $C$, there is a natural isomorphism $\mathcal{R}(B\bip C)\cong\mathcal{R}(B)\times\mathcal{R}(C)$.
Moreover, $\mathsf{B}(\mathcal{X}\times\mathcal{Y})\cong\mathsf{B}(\mathcal{X})\bip\mathsf{B}(\mathcal{Y})$
in $\mb$.
\item\label{it:itgrc4} The constructions from $B$ to $\mathcal{R}(B)$ and from $\mathcal{X}$ to $\mathsf{B}(\mathcal{X})$ are mutually inverse up to natural isomorphism, on objects and isomorphisms, and preserve products.
\end{enumerate}
\end{lemma}

\begin{proof}
(\ref{it:itgrc1}).
We first prove that the two equivalence relations defined in \eqref{eq:eqcp} are exactly $\rho_0$ and $\rho_1$.
Let $e=(a,b)\in E(B)$. By \eqref{eq:eqcer}, an edge $f=(c,d)$ belongs to $R_B[e]$ if and only if $(a,d)$ is an edge of $B$, or equivalently, $d\in N_B(a)$.
Thus $R_B[e]=\{(c,d)\in E(B):d\in N_B(a)\}$.
If $e\rho_0 f$, then the two edges have the same endpoint in $B^0$, so $R_B[e]=R_B[f]$. Conversely, suppose that $e=(a,b)$ and $f=(a',b')$ satisfy $R_B[e]=R_B[f]$.
By the preceding description of $R_B[e]$ and the fact that every vertex lies on an edge, we have $N_B(a)=N_B(a')$.
Since $B$ is twin-free, $a=a'$, and hence $e\rho_0f$.
Therefore, $\equiv_0$ is exactly $\rho_0$.
The same argument shows that $\equiv_1$ is exactly $\rho_1$.

We next verify that $\mathcal{R}(B)\in\mathfrak{G}$.
If $e\rho_0f$ and $e\rho_1f$, then $e$ and $f$ have both endpoints in common and therefore $e=f$. Hence $\rho_0\cap\rho_1$ is the diagonal relation, so $\mathcal{R}(B)$ is thin.

For every edge $e=(a,b)$ of $B$,  we have $e\rho_0e$ and $e\rho_1e$, so $eR_Be$.
Hence $R_B$ is reflexive. If two edges share an endpoint in either $B^0$ or $B^{1}$, then they  are related by $R_B$.
Since $B$ is connected, any two of its edges can be joined by a sequence of edges in which consecutive edges share an endpoint. Consecutive edges in this sequence are related by $R_B$, so $R_B$ is connected.

To verify rectangularity, let $e=(a,b)$ and $f=(c,d)$. By \eqref{eq:eqcer}, if $eR_Bf$, then  $g=(a,d)$ is an edge satisfying   $g\rho_0e$ and $g\rho_1f$. Conversely, if an edge $h$ satisfies $h\rho_0e$ and $h\rho_1f$, then $h=(a,d)$ and hence $eR_Bf$. Since $\rho_0$ and $\rho_1$ are the two equivalence relations from \eqref{eq:eqcp}, this is precisely rectangularity. Thus $\mathcal{R}(B)\in\mathfrak{G}$, which proves part~(\ref{it:itgrc1}).

(\ref{it:itgrc2}).
We first prove that $\mathcal R(\mathsf{B}(\mathcal{X}))\cong\mathcal{X}$.
Since $\mathcal{X}$ is thin, the construction of $\mathsf{B}(\mathcal{X})$ defines a bijection $\iota_X$ from $X$ to $E(\mathsf{B}(\mathcal{X}))$ with  $\iota_X(x)=e_x$.
In $\mathsf{B}(\mathcal{X})$, the edges $e_x$ and $e_z$ have the same endpoint in $(\mathsf{B}(\mathcal{X}))^0$ if and only if $z\equiv_0x$, and $e_z$ and $e_y$ have the same endpoint in $(\mathsf{B}(\mathcal{X}))^1$ if and only if $z\equiv_1y$. Therefore $e_xR_{\mathsf B(\mathcal{X})}e_y$ if and only if there is an element $z\in X$ such that $z\equiv_0x$ and $z\equiv_1y$.
By rectangularity, $e_xR_{\mathsf B(\mathcal{X})}e_y$ if and only if $xRy$. Thus $\iota_X$ is a relation isomorphism from $\mathcal{X}$ to $\mathcal{R}(\mathsf{B}(\mathcal{X}))$.

We now prove that $\mathsf{B}(\mathcal{X})$ is connected. Suppose first that $xRy$. By \eqref{eq:eqrectang}, there is an element $z\in X$ such that $e_x$ and $e_z$ share their endpoint in $\mathsf{B}(\mathcal{X})^0$, and $e_z$ and $e_y$ share their endpoint in $\mathsf{B}(\mathcal{X})^1$. Thus $e_x, e_z, e_y$, with consecutive repetitions removed, form a walk in the line graph of $\mathsf{B}(\mathcal{X})$. If $yRx$, the same argument produces a walk in the reverse direction.
Since $R$ is connected, any two elements $x, y$ of $X$ can be joined by a sequence $x=x_0, x_1,\ldots, x_m=y $ such that, for every $0\leq i\leq m-1$, either $x_iRx_{i+1}$ or $x_{i+1}Rx_{i}$. Replacing each step by the corresponding edge sequence gives a walk from $e_x$ to $e_y$ in the line graph of $\mathsf{B}(\mathcal{X})$.
Hence the line graph of $\mathsf{B}(\mathcal{X})$ is connected. Every vertex of  $\mathsf{B}(\mathcal{X})$ lies on an edge by construction, so $\mathsf{B}(\mathcal{X})$ is connected.

We next prove that $\mathsf{B}(\mathcal{X})$ is twin-free.
Let $([x]_{\equiv_0},0)$ be a vertex in $(\mathsf{B}(\mathcal{X}))^0$. By \eqref{eq:eqrectang},
$([y]_{\equiv_1},1)$ is adjacent to $([x]_{\equiv_0},0)$ if and only if there is an element $z\in X$ such that $z\equiv_0x$ and $z\equiv_1y$, which is equivalent to $xRy$. The set $R[x]$ is a union of $\equiv_1$ classes by the definition of $\equiv_1$. Thus
$N_{\mathsf{B}(\mathcal{X})}(([x]_{\equiv_0},0))=\{([y]_{\equiv_1},1):y\in R[x]\}$.
If the neighbourhoods of two vertices in $(\mathsf{B}(\mathcal{X}))^0$ are equal, then their $\equiv_0$ classes are equal. Similarly, $N_{\mathsf{B}(\mathcal{X})}(([x]_{\equiv_1},1))
 =\{([y]_{\equiv_0},0):y\in R^{-1}[x]\}$, and if the neighbourhoods of two vertices in $(\mathsf{B}(\mathcal{X}))^1$ are equal, then their $\equiv_1$ classes are equal.
Vertices in opposite colour classes also have different neighbourhoods, since these neighbourhoods are nonempty subsets of disjoint colour classes. Hence $\mathsf{B}(\mathcal{X})$ is twin-free.

Suppose that $f:\mathcal{X}\to\mathcal{Y}$ is a relation isomorphism. Then $xRy$ if and only if $f(x)Qf(y)$. Consequently, for $i=0,1$, the elements $x,x'\in X$ satisfy $x\equiv_i^R x'$ in $\mathcal{X}$ if and only if $f(x)$, $f(x')$ satisfy $f(x)\equiv_i^Q f(x')$ in $\mathcal{Y}$.
For each $i\in \{0,1\}$, define the induced graph map $\widehat{f}:\mathsf{B}(\mathcal{X})\to\mathsf{B}(\mathcal{Y})$ by
\[\widehat f(([x]_{\equiv_i^R},i))=([f(x)]_{\equiv_i^Q},i),\]
where the equivalence class on each side is taken in its respective relation structure. This map
is well defined and bijective. It maps the edge $e_x$ to the edge $e_{f(x)}$, and so it is an isomorphism in $\mb$. This completes part~(\ref{it:itgrc2}).

(\ref{it:itgrc3}).
 An edge of $B\bip C$ is uniquely specified by a pair
of edges, one of $B$ and one of $C$, with endpoints paired by colour.
Under this edge bijection, equality of a colour-$i$ endpoint is exactly
coordinatewise equality of colour-$i$ endpoints. Moreover, the two
witnesses for a composite relation can be chosen independently in the
factors. Hence $R_{B\bip C}$ corresponds to $R_B\times R_C$, giving
\begin{equation}\label{eq:eqrpp}
 \mathcal R(B\bip C)\cong\mathcal R(B)\times\mathcal R(C).
\end{equation}
By parts (\ref{it:itgrc1})--(\ref{it:itgrc2}) and
Lemmas~\ref{le:lebp}--\ref{le:lebpa}, this also shows that
$\mathcal X\times\mathcal Y$ belongs to $\mathfrak G$.

All row and column neighbourhoods of the reflexive factors are nonempty.
Projecting their Cartesian products gives, for $i=0,1$,
\[
 (x,y)\equiv_i^{R\times Q}(x',y')
 \iff x\equiv_i^Rx'\ \text{and}\ y\equiv_i^Qy'.
\]
The resulting map on vertices is
\[
 ([(x,y)]_{\equiv_i^{R\times Q}},i)
 \longmapsto
 \bigl(([x]_{\equiv_i^R},i),([y]_{\equiv_i^Q},i)\bigr).
\]
It is a bijection on each colour class and sends the edge indexed by
$(x,y)$ to the pair of edges indexed by $x$ and $y$. Thus it preserves
adjacency in both directions and proves
\begin{equation}\label{eq:eqbpp}
 \mathsf B(\mathcal X\times\mathcal Y)
 \cong\mathsf B(\mathcal X)\bip\mathsf B(\mathcal Y).
\end{equation}

(\ref{it:itgrc4})
For $a\in B^0$, choose an edge $e$ incident with $a$ and define $\vartheta_B(a)=([e]_{\rho_0},0)$.
Since all edges incident with $a$ belong to the same $\rho_0$ class, this definition is independent of the choice of $e$. Similarly, for $b\in B^1$, choose an incident edge $f$ and define $\vartheta_B(b)=([f]_{\rho_1},1)$.
This is likewise well defined. Part~(\ref{it:itgrc1}) shows that every $\rho_0$ class consists exactly of the edges incident with one vertex of $B^0$, and every $\rho_1$ class consists exactly of the edges incident with one vertex of $B^1$. Hence $\vartheta_B$ is a bijection that preserves the two colour classes.

If $e=(a,b)$ is an edge of $B$, then the edge of $\mathsf{B}(\mathcal{R}(B))$ corresponding to the element $e$ is $(([e]_{\rho_0},0),([e]_{\rho_1},1))=(\vartheta_B(a),\vartheta_B(b))$.
Conversely, every edge of $\mathsf{B}(\mathcal{R}(B))$ arises from a unique edge $e$ of $B$. Therefore $\vartheta_B$ is an isomorphism from $B$ to $\mathsf{B}(\mathcal{R}(B))$ in $\mb$.
Together with $\mathcal{X}\cong\mathcal{R}(\mathsf{B}(\mathcal{X}))$, Eqs.~\eqref{eq:eqrpp}, and \eqref{eq:eqbpp},  this shows that the two constructions are mutually inverse up to natural isomorphism and preserve products.
\end{proof}

\begin{lemma}\label{le:lefc}
Let $\mathcal{X}$ and $\mathcal{Y}$ be nonempty relation structures. If $\mathcal{X}\times\mathcal{Y}\in\mathfrak{G}$,
then $\mathcal{X},\mathcal{Y}\in\mathfrak{G}$.
\end{lemma}

\begin{proof}
Write $\mathcal{X}=(X,R)$ and $\mathcal{Y}=(Y,Q)$, where $X,Y$ are nonempty. Since $\mathcal{X}\times\mathcal{Y}$ is finite, both $X$ and $Y$ are finite. For $x\in X$ and $y\in Y$, reflexivity of $R\times Q$ gives $(x,y)(R\times Q)(x,y)$. Hence $xRx$ and $yQy$, so both factors are
reflexive.

To prove connectedness, take $x,y\in X$ and fix $z\in Y$. Since $R\times Q$ is connected, there is a sequence $(x,z)=(x_0,z_0), (x_1,z_1),\ldots, (x_m,z_m)=(y,z)$ such that, for every $0\leq i\leq m-1$, either $(x_i,z_i)(R\times Q)(x_{i+1},z_{i+1})$ or $(x_i,z_i)(R\times Q)^{-1}(x_{i+1},z_{i+1})$.
Projecting onto the first coordinate gives a sequence $x=x_0, x_1,\ldots, x_m=y$ such that, for every $0\leq i\leq m-1$, either $x_iRx_{i+1}$ or $x_iR^{-1}x_{i+1}$. Hence $R$ is connected. The same argument shows that $Q$ is connected.

Suppose that $\mathcal{X}$ is not thin. There are distinct $x,x'\in X$ for which $R[x]=R[x']$ and $R^{-1}[x]=R^{-1}[x']$. Fix $y\in Y$. Then
\[(R\times Q)[(x,y)]=R[x]\times Q[y]=R[x']\times Q[y]=(R\times Q)[(x',y)],\]
\[(R\times Q)^{-1}[(x,y)]=R^{-1}[x]\times Q^{-1}[y]=R^{-1}[x']\times Q^{-1}[y]=(R\times Q)^{-1}[(x',y)].\]
This contradicts the fact that $\mathcal{X}\times\mathcal{Y}$ is thin. Thus $\mathcal{X}$ is thin. A similar argument shows that $\mathcal{Y}$ is thin.

We prove rectangularity of $R$. Let $x,x'\in X$. If $xRx'$, fix $y\in Y$. Since $yQy$, we have $(x,y)(R\times Q)(x',y)$. Rectangularity of the product gives $(u,v)\in X\times Y$ such that
\[
\begin{aligned}
R[u]\times Q[v]&=R[x]\times Q[y],\\
R^{-1}[u]\times Q^{-1}[v]&=R^{-1}[x']\times Q^{-1}[y].
\end{aligned}
\]
All four factor neighbourhoods are nonempty by reflexivity, so projection of these Cartesian products yields $R[u]=R[x]$ and $R^{-1}[u]=R^{-1}[x']$. Hence $u\equiv_0^Rx$ and $u\equiv_1^Rx'$.

Conversely, suppose $u\in X$ satisfies $R[u]=R[x]$ and $R^{-1}[u]=R^{-1}[x']$. Fix $y\in Y$. The identities
\[
(R\times Q)[(u,y)]=(R\times Q)[(x,y)],\qquad
(R\times Q)^{-1}[(u,y)]=(R\times Q)^{-1}[(x',y)]
\]
show that $(u,y)\equiv_0^{R\times Q}(x,y)$ and $(u,y)\equiv_1^{R\times Q}(x',y)$ in the product relation. Its rectangularity gives $(x,y)(R\times Q)(x',y)$, hence $xRx'$. This proves rectangularity of $R$, and the same argument proves rectangularity of $Q$. Thus both factors belong to $\mathfrak G$.
\end{proof}

\begin{corollary}\label{co:cope}
Let $B$ be a connected twin-free ordered bipartite graph. Then $B$ is $f$-prime in $\mb$ if and only if $\mathcal{R}(B)$ is $d$-prime. Moreover, the one-element relation structure corresponds to $K_2$.
\end{corollary}

\begin{proof}
If $B\cong K_2$, then both $B$ and $\mathcal{R}(B)$ are units, and neither is prime. We may therefore assume that $B$ is nonunit.
Suppose first that $B$ has a nontrivial factorization $B\cong X\bip Y$
with neither factor isomorphic to $K_2$ in $\mb$.

The surjective coordinate projections show that both factors are connected, and Lemma~\ref{le:lebpa} shows that they are twin-free.
By Lemma~\ref{le:legrc}(\ref{it:itgrc3}), $\mathcal{R}(B)\cong\mathcal{R}(X)\times\mathcal{R}(Y)$.
 Since $\mathcal{R}(B)\in\mathfrak{G}$, Lemma \ref{le:lefc} implies that $\mathcal{R}(X),\mathcal{R}(Y)\in\mathfrak{G}$.
The relation structure associated with a connected ordered graph has one element if and only if the graph is $K_2$. Thus both relation factors are nonunit, and $\mathcal{R}(B)$ is not $d$-prime.

Conversely, suppose that $\mathcal{R}(B)\cong\mathcal{X}\times\mathcal{Y}$
with both $\mathcal{X}$ and $\mathcal{Y}$ nonunit. By Lemma~\ref{le:lefc},  $\mathcal{X},\mathcal{Y}\in\mathfrak{G}$. Applying Lemma~\ref{le:legrc}, we obtain
\[B\cong\mathsf{B}(\mathcal{R}(B))\cong\mathsf{B}(\mathcal{X}\times\mathcal{Y})\cong
 \mathsf{B}(\mathcal{X})\bip\mathsf{B}(\mathcal{Y}).
\]
Since $\mathcal{X}$ and $\mathcal{Y}$ are nonunit, neither graph factor is $K_2$. Thus $B$ is not $f$-prime. This proves the equivalence. The one-element relation structure corresponds to $K_2$ by construction.
\end{proof}

\section{Unique factorization for the fibre product over \texorpdfstring{$K_2$}{K2}}\label{sec:biproduct}
We now prove Theorems~\ref{th:thuf}, \ref{th:thci}, and \ref{th:thfc}, beginning with unique factorization for relation structures.
\subsection{Unique factorization of relation structures}\label{subsec:factorization}
The following result is a formulation of \cite[Theorem~4.1]{McK71} in terms of coordinate isomorphisms. For completeness, we state this formulation and derive it from the cited result.

\begin{lemma}\label{th:thsr}
Let $\mathcal{A}$ be a finite connected reflexive thin relation structure. Suppose $r,s\geq 1$ and
\[
 \alpha:\mathcal{A}\longrightarrow\prod_{i=1}^r\mathcal{A}_i,
 \qquad
 \beta:\mathcal{A}\longrightarrow\prod_{j=1}^s\mathcal{C}_j
\]
are relation isomorphisms. Then there exist relation structures $\mathcal D_{ij}$ (possibly units), relation isomorphisms
\[
 \lambda_i:\mathcal{A}_i\longrightarrow\prod_{j=1}^s\mathcal D_{ij},
 \qquad
 \mu_j:\mathcal{C}_j\longrightarrow\prod_{i=1}^r\mathcal D_{ij},
\]
and the canonical transposition of coordinates
\[
 \tau:\prod_{j=1}^s\prod_{i=1}^r\mathcal D_{ij}
 \longrightarrow
 \prod_{i=1}^r\prod_{j=1}^s\mathcal D_{ij}
\]
such that
\begin{equation}\label{eq:eqsc}
 \left(\prod_{i=1}^r\lambda_i\right)\circ\alpha
 =
 \tau\circ\left(\prod_{j=1}^s\mu_j\right)\circ\beta.
\end{equation}
\end{lemma}

\begin{proof}
Let $\mathcal{A}=(A, R)$. For $1\leq i\leq r$ and $1\leq j\leq s$, let $p_i$ denote the $i$-th projection from $\prod_{t=1}^r\mathcal{A}_t$ to $\mathcal{A}_i$, and let $q_j$ denote the $j$-th projection from $\prod_{t=1}^s\mathcal{C}_t$ to $\mathcal{C}_j$. For $1\leq i\leq r$ and $1\leq j\leq s$, we define equivalence relations $E_i$ and $F_j$ on $A$ by
\begin{align*}
 aE_i b
 &\quad\Longleftrightarrow\quad
  p_i(\alpha(a))= p_i(\alpha(b)),\\
 aF_j b
 &\quad\Longleftrightarrow\quad
  q_j(\beta(a))= q_j(\beta(b)).
\end{align*}
For each $1\leq i\leq r$, define the quotient structure  $\mathcal{A}/E_i$ on the set of $E_i$ classes by
\[[a]_{E_i}\,\overline R_i\,[b]_{E_i}\quad\Longleftrightarrow\quad p_i(\alpha(a))\,R_i\, p_i(\alpha(b)),
\]
where $R_i$ is the relation of $\mathcal{A}_i$.
 For every $1\leq i\leq r$, define a map $\overline{\alpha}_i$ from $\mathcal{A}/E_i$ to $\mathcal{A}_i$ by
 $\overline{\alpha}_i([a]_{E_i})= p_i(\alpha(a))$
 for $[a]_{E_i}\in \mathcal{A}/E_i$.
 Since $\alpha$ and $ p_i$ are surjective, the map $\overline{\alpha}_i$ is surjective.
For any elements $[a]_{E_i}, [b]_{E_i}\in \mathcal{A}/E_i$, $\overline{\alpha}_i([a]_{E_i})=\overline{\alpha}_i([b]_{E_i})$ if and only if $ p_i(\alpha(a))= p_i(\alpha(b))$.
 The latter equality holds if and only if $aE_i b$, equivalently, $[a]_{E_i}=[b]_{E_i}$. Hence $\overline{\alpha}_i$ is injective and therefore bijective. The definition of $\overline{R}_i$ shows that $\overline{\alpha}_i$ is a relation isomorphism. Similarly, we construct a relation isomorphism $\overline{\beta}_j$ from $\mathcal{A}/F_j$ to $\mathcal{C}_j$ such that $\overline{\beta}_j([a]_{F_j})= q_j(\beta(a))$ for $1\leq j\leq s$.

Since $R$ is reflexive, both $R\circ R^{-1}$ and $R^{-1}\circ R$ contain $R\cup R^{-1}$. Thus both compositions are connected whenever $R$ is connected.
By Definition~\ref{de:derds}, $\mathcal{A}$ is thin. Therefore McKenzie's theorem~\cite[Theorem~4.1]{McK71}
applies and gives the strict refinement property for factor relations.

Applying this property to  $(E_i)_{1\leq i\leq r}$ and $(F_j)_{1\leq j\leq s}$, we obtain  an $r\times s$ array of factor relations $E_{ij}$ with $1\leq i\leq r$ and $1\leq j\leq s$ such that
\begin{equation*}
 E_i=\cap_{j=1}^s E_{ij},
 \qquad
 F_j=\cap_{i=1}^r E_{ij},
\end{equation*}
and the natural maps
\begin{align*}
 \xi_i:\mathcal{A}/E_i&\longrightarrow\prod_{j=1}^s\mathcal{A}/E_{ij},&[a]_{E_i}
 &\longmapsto([a]_{E_{i1}},\ldots,[a]_{E_{is}}),\\
\zeta_j:\mathcal{A}/F_j&\longrightarrow\prod_{i=1}^r\mathcal{A}/E_{ij},&[a]_{F_j}
 &\longmapsto ([a]_{E_{1j}},\ldots,[a]_{E_{rj}})
\end{align*}
are relation isomorphisms for $1\leq i\leq r$ and $1\leq j\leq s$.

For each pair $(i,j)$, let $\mathcal{D}_{ij}=\mathcal{A}/E_{ij}$.
Define the maps
\[\lambda_i =\xi_i\circ\overline\alpha_i^{-1}
 :\mathcal{A}_i\longrightarrow\prod_{j=1}^s\mathcal{D}_{ij}
\]
and
\[\mu_j =\zeta_j\circ\overline\beta_j^{-1}
 :\mathcal{C}_j\longrightarrow\prod_{i=1}^r\mathcal{D}_{ij}.
\]
Since $\xi_i$ and $\overline{\alpha}^{-1}_i$ are relation isomorphisms, $\lambda_i$ is also a relation isomorphism. Similarly, $\mu_j$ is a relation isomorphism. For every $a\in A$,  these maps satisfy
 \[\lambda_i\bigl( p_i(\alpha(a))\bigr)
 =\bigl([a]_{E_{i1}},\ldots,[a]_{E_{is}}\bigr),\quad
 \mu_j\bigl( q_j(\beta(a))\bigr)
 =\bigl([a]_{E_{1j}},\ldots,[a]_{E_{rj}}\bigr).\]
The maps $\lambda_i$ and $\mu_j$ are relation isomorphisms
as compositions of relation isomorphisms. For every $a\in A$, the $j$-th coordinate of $\lambda_i( p_i(\alpha(a)))$ and the $i$-th coordinate of $\mu_j( q_j(\beta(a)))$ are both $[a]_{E_{ij}}$. After the transposition $\tau$, the two sides of Eq.~\eqref{eq:eqsc} therefore agree in every coordinate. This proves the asserted compatibility for the given maps $\alpha$ and $\beta$.
\end{proof}

\begin{remark}\label{re:reunit}
In Lemma~\ref{th:thsr}, a universal factor equivalence $E_{ij}$ gives the one-element reflexive quotient. Such unit entries can be omitted from the displayed row and column products.
\end{remark}

Let $\mathcal{A}$ be a finite connected reflexive thin relation structure.
For a rectangular array $(\mathcal{D}_{ij})$ obtained in Lemma~\ref{th:thsr}, the family $(\mathcal{D}_{ij})_{1\leq j\leq s}$ with $i$ fixed is the $i$-th \emph{row}, and $\prod_{j=1}^s\mathcal{D}_{ij}$ is its \emph{row product}. Similarly, the family $(\mathcal{D}_{ij})_{1\leq i\leq r}$ with $j$ fixed is the $j$-th \emph{column}, and $\prod_{i=1}^r\mathcal{D}_{ij}$
is its \emph{column product}.

\subsection{Prime factorization of ordered bipartite graphs}\label{subsec-3.2}
We now apply the preceding results to prove Theorem~\ref{th:thuf}.
\begin{proof}[Proof of Theorem~\ref{th:thuf}]
If $B\cong K_2$, every displayed factor is connected by surjectivity of its projection. Since the number of edges of a fibre product is the product of the factor edge counts, each displayed factor has one edge and is $K_2$. For part~(\ref{it:ituf1}), take $\mathcal S=\varnothing$ and all indexed
parts empty. For part~(\ref{it:ituf2}), a decomposition into nonunit factors must be empty, so $r=s=0$. We may assume $B$ nonunit, which forces $r,s\geq1$.

(\ref{it:ituf1}).
Fix the two factorization isomorphisms.
For each $1\leq i\leq r$ and $1\leq j\leq s$, let $ p_i$ and $ q_j$ be the coordinate projections from $B$ to $G_i$ and $H_j$, respectively. Since $ p_i$ and $ q_j$ are surjective, all $G_i$ and $H_j$ are connected.
Repeated application of Lemma~\ref{le:lebpa} shows that every $G_i$ and $H_j$ is twin-free. Let $\mathcal{A}=\mathcal{R}(B)$.
By Lemma~\ref{le:legrc}(\ref{it:itgrc3}),
$\mathcal{R}(B)\cong \mathcal{R}(G_1)\times\cdots\times \mathcal{R}(G_r)\cong \mathcal{R}(H_1)\times\cdots\times \mathcal{R}(H_s)$. Choose isomorphisms
\[
 \alpha:\mathcal{A}\longrightarrow\prod_{i=1}^r\mathcal{R}(G_i),
 \qquad
 \beta:\mathcal{A}\longrightarrow\prod_{j=1}^s\mathcal{R}(H_j).
\]
Applying Lemma~\ref{th:thsr}, we obtain an $r\times s$ array $(\mathcal{D}_{ij})$ such that
\begin{equation}\label{eq:eqcor}
 \mathcal{R}(G_i)\cong\prod_{j=1}^s\mathcal{D}_{ij},
 \qquad
 \mathcal{R}(H_j)\cong\prod_{i=1}^r\mathcal{D}_{ij}.
\end{equation}
Moreover, $\mathcal{A}\cong\prod_{i=1}^r\prod_{j=1}^s\mathcal{D}_{ij}\cong\prod_{j=1}^s\prod_{i=1}^r\mathcal{D}_{ij}$.
Since the $i$-th row product in ~\eqref{eq:eqcor} is isomorphic to $\mathcal{R}(G_i)\in\mathfrak{G}$, repeated application of Lemma~\ref{le:lefc} shows that each entry  $\mathcal{D}_{ij}$ belongs to $\mathfrak{G}$.

We next factor each nonunit $\mathcal{D}_{ij}$ into $d$-prime factors. In a nontrivial direct product factorization of a finite relation structure, each nonunit factor has strictly fewer elements than the relation structure.
Repeated factorization must terminate and yields a factorization of each row product $\mathcal{R}(G_i)$ into $d$-prime factors.
Repeated application of Lemma~\ref{le:lefc} ensures that every factor belongs to $\mathfrak{G}$.

Index the $d$-prime nonunit factors of all $\mathcal{D}_{ij}$, for $1\leq i\leq r, 1\leq j\leq s$, as $\mathcal T_\lambda$ with $\lambda\in \mathcal S$.
Each indexed occurrence $\mathcal T_\lambda$ arises from a unique entry $\mathcal{D}_{ij}$.
Let $U_i$ be the set of indices arising from row $i$, and define $V_j$ analogously for column $j$. Since every $\mathcal{D}_{ij}$ belongs to one row and one column, the indexed families $\{U_i:1\leq i\leq r\}$ and $\{V_j:1\leq j\leq s\}$ both partition $\mathcal S$, with empty parts allowed.

Define $B_\lambda=\mathsf{B}(\mathcal T_\lambda)$.
Lemma~\ref{le:legrc} shows that each $B_\lambda$ is connected and twin-free. Corollary~\ref{co:cope} shows that each is $f$-prime. Applying $\mathsf{B}$ to the row products, column products, and the product of all $\mathcal{D}_{ij}$, and using Lemma~\ref{le:legrc}(\ref{it:itgrc3}) and (\ref{it:itgrc4}), we obtain
\[
  G_i\cong\Bipprod_{\lambda\in U_i}B_\lambda,
 \qquad
 H_j\cong\Bipprod_{\lambda\in V_j}B_\lambda,
 \qquad
 B\cong\Bipprod_{\lambda\in \mathcal S}B_\lambda.
\]
This proves part~(\ref{it:ituf1}).

(\ref{it:ituf2}).
Assume  that all $G_i$ and $H_j$ are $f$-prime for $1\leq i\leq r$ and $1\leq j\leq s$. Then every row product $\mathcal{R}(G_i)$ in \eqref{eq:eqcor} is $d$-prime and nonunit by Corollary~\ref{co:cope}. Therefore its row contains  exactly one nonunit entry (since two nonunit entries would give a nontrivial factorization of the row product). By the same argument, every column contains exactly one nonunit entry.

Consequently, matching each row to the column containing its nonunit entry defines a bijection between the rows and columns. Thus $r=s$. Let $\sigma(i)$ be the unique column containing the nonunit entry in row $i$. By the row and column isomorphisms in \eqref{eq:eqcor}, we have
\[ \mathcal{R}(G_i)\cong\mathcal{D}_{i,\sigma(i)}
 \cong\mathcal{R}(H_{\sigma(i)}),\]
where $\mathcal{D}_{i,\sigma(i)}$ is the unique nonunit entry in row $i$, corresponding to $j=\sigma(i)$.
By Lemma~\ref{le:legrc}, $G_i\cong H_{\sigma(i)}$ in $\mb$. This proves part~(\ref{it:ituf2}) and the uniqueness statement.
\end{proof}


\subsection{Coordinate isomorphisms and colour signs}\label{subsec-3.3}
We next prove Theorems~\ref{th:thci} and \ref{th:thfc}.

\begin{proof}[Proof of Theorem~\ref{th:thci}]
Since $B_i$ and $C_j$, for $1\leq i\leq r$ and $1\leq j\leq s$, are nonunit factors, both $B$ and $C$ are nonunit.
By Lemma~\ref{le:legrc}(\ref{it:itgrc3}), we have the natural isomorphisms
\[
 \eta_B:\mathcal{R}(B)\longrightarrow\prod_{i=1}^r\mathcal{R}(B_i),
 \qquad
 \eta_C:\mathcal{R}(C)\longrightarrow\prod_{j=1}^s\mathcal{R}(C_j).
\]
The graph isomorphism $\varphi$ maps an edge $e=\{{\bf x},{\bf y}\}$ of $B$ to the edge $\varphi(e)=\{\varphi({\bf x}),\varphi({\bf y})\}$ of $C$. Since $\varphi$ preserves the two colour classes, it maps each endpoint in $B^0$ to an endpoint in $C^0$ and each endpoint in $B^1$ to an endpoint in $C^1$. Hence it induces a relation isomorphism $\Phi$ from $\mathcal{R}(B)$ to $\mathcal{R}(C)$.

 Applying Lemma~\ref{th:thsr} to the two coordinatizations $\alpha=\eta_B$ and $\beta=\eta_C\circ\Phi$ of $\mathcal{R}(B)$, we obtain relation structures $\mathcal{D}_{ij}$ satisfying
\[
 \mathcal{R}(B_i)\cong\prod_{j=1}^s\mathcal{D}_{ij},
 \qquad
 \mathcal{R}(C_j)\cong\prod_{i=1}^r\mathcal{D}_{ij}.\]
By Corollary~\ref{co:cope}, every row product $\prod_{j=1}^s\mathcal{D}_{ij}$ and every column product $\prod_{i=1}^r\mathcal{D}_{ij}$ is $d$-prime and nonunit. As shown in the proof of Theorem~\ref{th:thuf}, each row and each column contains exactly one nonunit entry. It follows that $r=s$, and there is a permutation $\sigma\in S_r$ such that the unique nonunit entry in row $i$ is $\mathcal{D}_{i,\sigma(i)}$.

After deleting all unit coordinates, the row and column isomorphisms in Lemma~\ref{th:thsr} become
\[
 \lambda_i:\mathcal{R}(B_i)\longrightarrow\mathcal{D}_{i,\sigma(i)},
 \qquad
 \mu_{\sigma(i)}:\mathcal{R}(C_{\sigma(i)})
 \longrightarrow\mathcal{D}_{i,\sigma(i)}.
\]
Define
\[
 \Phi_i=\mu_{\sigma(i)}^{-1}\circ\lambda_i:
 \mathcal{R}(B_i)\longrightarrow\mathcal{R}(C_{\sigma(i)}).
\]

Let $p_i:\prod_{h=1}^r\mathcal{R}(B_h)\to\mathcal{R}(B_i)$ and $q_j:\prod_{h=1}^r\mathcal{R}(C_h)\to\mathcal{R}(C_j)$ denote the source and target coordinate projections, respectively.
Combining this with Eq.~\eqref{eq:eqsc}, we have
\begin{equation}\label{eq:eqeci}
 q_{\sigma(i)}\circ\eta_C\circ\Phi=\Phi_i\circ p_i\circ\eta_B.
\end{equation}
Thus the $\sigma(i)$-th edge coordinate of $\Phi(e)$
is obtained by applying the fixed relation isomorphism $\Phi_i$ to
 the $i$-th edge coordinate of $e$.

By Lemma~\ref{le:legrc}(\ref{it:itgrc2}), the relation isomorphism $\Phi_i$ induces a unique graph isomorphism $\varphi_i$ from $B_i$ to $C_{\sigma(i)}$, which preserves the ordered bipartition. We now show that these factor maps also determine the vertex coordinates of $\varphi$.

Fix a vertex ${\bf b}=(b_1,\ldots,b_r)\in V(B)$. Since  $B$ is connected and has an edge, ${\bf b}$ has a neighbour. Choose an incident edge ${\bf e}=(e_1,\ldots,e_r)$,
where $e_i$ is an edge of $B_i$ incident with $b_i$. The edge $\Phi(\mathbf e)$ is incident with $\varphi(\mathbf b)$. By Eq.~\eqref{eq:eqeci}, the $\sigma(i)$-th edge coordinate of $\Phi({\bf e})$ is $\Phi_i(e_i)$. The graph isomorphism induced by $\Phi_i$ sends the endpoint $b_i$ of $e_i$ to the endpoint $\varphi_i(b_i)$ of $\Phi_i(e_i)$ having the same colour. Since $\varphi$ preserves the global colour, these endpoints are the coordinates of $\varphi({\bf b})$. Therefore, if $\varphi({\bf b})=(c_1,\ldots,c_r)$, then $c_{\sigma(i)}=\varphi_i(b_i)$
for every $i$, proving Eq.~\eqref{eq:eqcis}. This completes the proof.
\end{proof}



\begin{proof}[Proof of Theorem~\ref{th:thfc}]
Since $B$ and $C$ are connected bipartite graphs, their bipartitions are unique up to interchange.  Thus there is a unique $\delta\in\Ftwo$ such that for every ${\bf b}\in V(B)$
\begin{equation}\label{eq:eqgsp}
 \mrc_C(\varphi({\bf b}))=\mrc_B({\bf b})+\delta.
\end{equation}
Both $B$ and $C$ are nonunit by the assumed nonunit factorizations. We use the fixed factorizations to identify their vertices with compatible coordinate tuples.
If $\delta=0$, then by Eq.~\eqref{eq:eqgsp}, $\varphi$ is an isomorphism from $B$ to $C$ in $\mb$. By Theorem~\ref{th:thci}, $r=s$ and there are a permutation $\sigma$ and factor isomorphisms that preserve colours and satisfy the coordinate formula. These maps satisfy Eqs.~\eqref{eq:eqfas} and \eqref{eq:eqocis} with $\delta=0$.

Now assume $\delta=1$. Regard $\varphi$ as a map into $C^{\op}$. Its colour map is $\mrc_{C^{\op}}=\mrc_C+1$, so Eq.~\eqref{eq:eqgsp} becomes $\mrc_{C^{\op}}(\varphi(b))=\mrc_B(b)$.
Thus $\varphi:B\to C^{\op}$ is an isomorphism in $\mb$. Since reversing every factor colouring reverses the global colouring, we have the following isomorphism in $\mb$,
\begin{equation}\label{eq:eqraf}
 C^{\op}
 \cong
 C_1^{\op}\bip\cdots\bip C_s^{\op}
 \quad\text{in }\mb.
\end{equation}
Reversing all factor colours in a nontrivial factorization of $C_j^{\op}$ produces a nontrivial factorization of $C_j$. Since $C_j$ is $f$-prime, each $C_j^{\op}$ is also $f$-prime for $1\leq j\leq s$.

 By Theorem~\ref{th:thci} and the decomposition \eqref{eq:eqraf}, we obtain $r=s$, a permutation $\sigma\in S_r$, and isomorphisms $\varphi_i: B_i\longrightarrow C_{\sigma(i)}^{\op}$ that preserve colours.
Viewed as ordinary isomorphisms into $C_{\sigma(i)}$, these maps interchange the two colour classes and hence satisfy Eq.~\eqref{eq:eqfas} with $\delta=1$. The coordinate formula~\eqref{eq:eqocis} follows from Theorem~\ref{th:thci}.

The sign $\delta$ is unique since it is determined by the ordinary isomorphism $\varphi$ in Eq.~\eqref{eq:eqgsp}. This completes the proof.
\end{proof}

The following example illustrates Theorem \ref{th:thfc}.
\begin{example}
Let $P$ and $Q$ be nonisomorphic connected twin-free nonunit $f$-prime ordered bipartite graphs. Let $B=P\bip Q$ and $C=Q\bip P$.
The coordinate transposition $(p,q)\longmapsto(q,p)$
is an isomorphism in $\mb$. The permutation in Theorem~\ref{th:thfc} interchanges the two positions (the factor $P$ of $B$ corresponds to the second factor of $C$, and the factor $Q$ corresponds to the first).

The identity map on vertices is also an ordinary isomorphism $\Und(B)\longrightarrow\Und(B^{\op})$
with global sign $1$. Since $B^{\op}\cong P^{\op}\bip Q^{\op}$,
both factor maps reverse their bipartitions.
 If one factor map preserves colours and another reverses them, the resulting coordinate map does not preserve the condition that all coordinates have the same colour. This explains why all factor maps must have the same sign in Theorem~\ref{th:thfc}.
\end{example}

We conclude this section with an application of Theorem \ref{th:thfc}.

\begin{corollary}\label{cor:corsfc}
Let $B_1,\ldots,B_r$ and $C_1,\ldots,C_s$ be connected twin-free nonunit $f$-prime ordered bipartite graphs, where $r, s\geq 1$. For $\boldsymbol{\epsilon}\in\Ftwo^r$ and $\boldsymbol{\eta}\in\Ftwo^s$, write $\Omega_{\ole}^B=\Und\bigl(B_1^{(\epsilon_1)}\bip\cdots\bip B_r^{(\epsilon_r)}\bigr)$,
and
$\Omega_{\overline{\boldsymbol{\eta}}}^C=\Und\bigl(C_1^{(\eta_1)}\bip\cdots\bip C_s^{(\eta_s)}\bigr)$.
Then $\Omega_{\ole}^B\cong\Omega_{\ol{\boldsymbol{\eta}}}^C$ if and only if $r=s$ and there are $\delta\in\Ftwo$ and a permutation $\sigma\in S_r$ such that, for each $1\leq i\leq r$, the following holds in $\mb$,
\begin{equation}\label{eq:eqsfc}
 B_i^{(\epsilon_i)}\cong C_{\sigma(i)}^{(\eta_{\sigma(i)}+\delta)}.
\end{equation}
\end{corollary}

\begin{proof}
First suppose that $\Omega_{\ole}^B$ is isomorphic to $\Omega_{\ol{\boldsymbol{\eta}}}^C$.
Equip $\Omega_{\ole}^B$ and $\Omega_{\ol{\boldsymbol{\eta}}}^C$ with the ordered bipartitions of $B_1^{(\epsilon_1)}\bip\cdots\bip B_r^{(\epsilon_r)}$ and $C_1^{(\eta_1)}\bip\cdots\bip C_s^{(\eta_s)}$, respectively.
By Theorem~\ref{th:thfc}, $r=s$, and there are a coordinate permutation $\sigma$ and a common global sign $\delta$.
By Eq.~\eqref{eq:eqfas}, each factor map from $B_i^{(\epsilon_i)}$ to $C_{\sigma(i)}^{(\eta_{\sigma(i)})}$ with sign $\delta$ preserves the ordered bipartition when its codomain is given the colouring of $C_{\sigma(i)}^{(\eta_{\sigma(i)}+\delta)}$. Thus $B_i^{(\epsilon_i)}\cong C_{\sigma(i)}^{(\eta_{\sigma(i)}+\delta)}$ in $\mb$, which proves Eq.~\eqref{eq:eqsfc}.

Conversely, take the fibre product of the factor isomorphisms in Eq.~\eqref{eq:eqsfc} and permute the output  coordinates by $\sigma$. The resulting map is an isomorphism in $\mb$ from $B_{1}^{(\epsilon_{1})}\bip\cdots\bip B_{r}^{(\epsilon_{r})}$ to $C_{1}^{(\eta_{1}+\delta)}\bip\cdots\bip C_{r}^{(\eta_{r}+\delta)}$. Forgetting the order of the bipartition gives an ordinary isomorphism from $\Omega_{\ole}^B$ to  $\Omega_{\ol{\boldsymbol{\eta}}}^C$.
This completes the proof.
\end{proof}

\section{Automorphism groups of the direct product of bipartite graphs}\label{subsec:automorphisms}
We first prove the component and group statements in Theorem~\ref{th:thrfa}.
We then specialize to distinct factors and to fibre powers, making the orbit intersections explicit in the latter case. The final applications are a stability transfer and a recovery of known double cover formulas.

\subsection{Proof of Theorem~\ref{th:thrfa}} We prove the component identification, orbit criterion, stabilizer formula, and group decomposition in that order. The distinction between an individual component and a coordinate factor will be used throughout.

\begin{proof}[Proof of Theorem~\ref{th:thrfa}]
All decompositions $\theta_i$ are fixed as in~\eqref{eq:eqcr}. We use ${\bf x},{\bf y}, {\bf z}$ for vertices of $P$ and ${\bf u},{\bf v}, {\bf w}$ for vertices of $\widetilde{P}$. Binary vectors $\mathbf{a}$, $\mathbf{b}$   and their
classes $\overline{\mathbf a}$ and $\overline{\mathbf b}$ label components, not vertices.

(\ref{it:ierfaa}).
We first prove that $\Theta$ is an isomorphism from $P$ to $\sqcup_{\ol{{\bf a}}\in \mathcal L_{\mathcal I}}\Lambda_{\ol{{\bf a}}}$. Let $\mathbf{x}=(x_1,\ldots,x_r)\in V(P)$, and write
$\Theta(\mathbf{x})=\mathbf{u}=(u_1,\ldots,u_N)$.  By the definition
of $\Theta$, for each $1\leq i\leq r$, $\theta_i(x_i)=(u_k)_{k\in I_i}$.
Since each $\theta_i$ is an isomorphism in $\mb$, it preserves the
ordered colour classes.  Consequently, for every $i\in\{1,\ldots,r\}$ and $k\in I_i$, $\mrc_{C_k}(u_k)=\mrc_{B_i}(x_i)$.
Consider the following set of vertices of the refined direct product:
\[
W=\bigl\{\mathbf v=(v_1,\ldots,v_N)\in V(\widetilde P):
 \mrc_{C_k}(v_k)=\mrc_{C_\ell}(v_\ell)
 \text{ whenever }k,\ell\in I_i\text{ for some }i\bigr\}.\]

The colour-preservation identity above shows that $\Theta(V(P))\subseteq W$.
Conversely,  take $\mathbf{v}\in W$.  For every $i$, the coordinates
$(v_k)_{k\in I_i}$ have a common colour and hence form a vertex of
$\Bipprod_{k\in I_i}C_k$.  Since $\theta_i$ is bijective, there is  a
unique $y_i\in V(B_i)$ with $\theta_i(y_i)=(v_k)_{k\in I_i}$.  The resulting tuple
$\mathbf{y}=(y_1,\ldots,y_r)$ belongs to $V(P)$ and satisfies
$\Theta(\mathbf{y})=\mathbf{v}$.  Uniqueness of each $y_i$ makes this
preimage unique.  Thus $\Theta$ is a bijection  from $V(P)$ onto $W$.

To check adjacency, let $\mathbf{z},\mathbf{z'}\in V(P)$, and put
$\mathbf{w}=\Theta(\mathbf{z})$ and
$\mathbf{w'}=\Theta(\mathbf{z'})$.  Since each $\theta_i$ is a graph isomorphism, we have
\[\begin{aligned}
 \mathbf{z}\sim_P\mathbf{z'}
 &\iff z_i\sim_{\Und(B_i)}z_i'\quad\text{for all }i\\
 &\iff (w_k)_{k\in I_i}
     \sim_{\Bipprod_{k\in I_i}C_k}(w'_k)_{k\in I_i}
     \quad\text{for all }i\\
 &\iff w_k\sim_{\Und(C_k)}w'_k\quad\text{for all }k\\
 &\iff \mathbf{w}\sim_{\widetilde P}\mathbf{w'}.
\end{aligned}
\]
Hence $\Theta$ is a graph isomorphism from $P$ onto the subgraph
of $\widetilde{P}$ induced by $W$.

We next identify this induced subgraph as a union of whole connected
components.  For any $\mathbf{v}=(v_1,\ldots,v_N)\in V(\widetilde{P})$,
define its coordinate colour vector by
\[
\widetilde{\pi}(\mathbf{v})
   =\bigl(\mrc_{C_1}(v_1),\ldots,\mrc_{C_N}(v_N)\bigr)
   \in\Ftwo^N.
\]
By Proposition~\ref{pro:procd}, the connected component containing $\mathbf{v}$ has label $\overline{\widetilde{\pi}(\mathbf{v})}\in Q_N$. For $N=1$ the same assertion is immediate.  By the definition of
$\mathcal L_{\mathcal I}$, this label belongs to $\mathcal L_{\mathcal I}$ if and only if $\widetilde{\pi}(\mathbf{v})$ is constant on each block $I_i$. Indeed, adding $\mathbf{1}_N$ changes every entry by the same amount and therefore does not affect the blockwise equalities.  It follows
that $W=\sqcup_{\ol{{\bf a}}\in \mathcal L_{\mathcal I}}V(\Lambda_{\ol{{\bf a}}})$.

Each $\Lambda_{\ol{{\bf a}}}$ is a connected component of $\widetilde{P}$,
so there are no edges between different such components.  Then we have
\[\Theta:P\longrightarrow
 \sqcup_{\ol{{\bf a}}\in \mathcal L_{\mathcal I}}\Lambda_{\ol{\bf a}}.\]
It remains to identify the restriction to each original component.
Recall that the coordinate colour vector of
$\mathbf{x}=(x_1,\ldots,x_r)\in V(P)$ is
$\pi(\mathbf{x})=(\mrc_{B_1}(x_1),\ldots,\mrc_{B_r}(x_r))$.
The colour-preservation identity gives $\widetilde{\pi}(\Theta(\mathbf{x}))=d(\pi(\mathbf{x}))$.
We also verify that the induced repetition map
$\ol{d}:Q_r\to Q_N$ is injective.  Suppose
$\ol{d}(\ole)=\ol{d}(\ol{\boldsymbol{\eta}})$.  Then
$d(\boldsymbol{\epsilon}+\boldsymbol{\eta})=t\mathbf{1}_N$ for some $t\in\Ftwo$.
Since every block $I_i$ is nonempty, comparing any coordinate in
$I_i$ gives $\epsilon_i+\eta_i=t$ for each $i$.  Hence
$\boldsymbol{\epsilon}+\boldsymbol{\eta}=t\mathbf1_r$, and $\ole=\ol{\boldsymbol{\eta}}$, as required.

Fix $\ole\in Q_r$.  Proposition~\ref{pro:procd}
(and the evident case $r=1$) shows that
$\mathbf{y}\in V(\Omega_{\ole})$ precisely when
$\overline{\pi(\mathbf{y})}=\ole$.  Taking quotient classes in this identity, we see that $\Theta$ maps
$\Omega_{\ole}$ into
$\Lambda_{\ol{d}(\ole)}$.
Conversely,  let $\mathbf{v}\in V(\Lambda_{\ol{d}(\ole)})$.
The vertex $\mathbf{v}$ lies in $W$ and hence has a unique preimage $\mathbf{y}\in V(P)$.  Thus
$\ol{d}(\overline{\pi(\mathbf{y})})=\ol{d}(\ole)$.
Since  $\ol{d}$ is injective, $\overline{\pi(\mathbf{y})}=\ole$, so
$\mathbf{y}\in V(\Omega_{\ole})$.
Thus \[
\Theta_{\ole}:\Omega_{\ole}
    \longrightarrow\Lambda_{\ol{d}(\ole)},
\] is a graph isomorphism.
Finally, since $\mathcal L_{\mathcal I}=\ol d(Q_r)$ and $\ol d$ is injective,
\[
|\mathcal L_{\mathcal I}|=|Q_r|=2^{r-1}.
\]
This proves part~(\ref{it:ierfaa}).

(\ref{it:ierfab}).
We first identify  the action of a factor symmetry on a refined component. Let $w=w_{\sigma,\mathbf f}\in\widetilde{\mfw}$, where
$\sigma\in S_N$ and
$f_k:\Und(C_k)\to\Und(C_{\sigma(k)})$
is an isomorphism for every $1\leq k\leq N$.
Since the factors are connected and bipartite, each $f_k$ has a unique colour sign $s_k\in\Ftwo$ satisfying, for any $z\in V(C_k)$,
\begin{equation}\label{eq:eqth16ps}
 \mrc_{C_{\sigma(k)}}(f_k(z))
   =\mrc_{C_k}(z)+s_k.
\end{equation}
Since $\mrc_{C_{\sigma(k)}}(f_k(z))+\mrc_{C_k}(z)$ is unchanged across every
edge and $C_k$ is connected, the sign $s_k$ is independent of $z$. Write $\mathbf{s}=(s_1,\ldots,s_N)\in\Ftwo^N$.

Fix $\mathbf{a}\in\Ftwo^N$ and a vertex
$\mathbf{u}=(u_1,\ldots,u_N)\in V(\Lambda_{\ol{\mathbf{a}}})$.
There is
$t\in\Ftwo$ such that
$\mrc_{C_k}(x_k)=a_k+t$ for every $k$.
Write $\mathbf{v}=w(\mathbf u)=(v_1,\ldots,v_N)$.
Then $v_{\sigma(k)}=f_k(x_k)$ and Eq.~\eqref{eq:eqth16ps} gives \[\widetilde\pi(\mathbf v)=\sigma_*(\mathbf a+\mathbf s)+t\mathbf1_N.\]
 Applying the same argument to $w_{\sigma,\mathbf{f}}^{-1}$, we have
\begin{equation}\label{eq:eqth16pca}
 w_{\sigma,\mathbf{f}}(\Lambda_{\ol{{\bf a}}})=\Lambda_{\overline{\sigma_*({\bf a}+{\bf s})}}.
\end{equation}
%

Let
$\varphi:\Lambda_{\ol{\mathbf a}}
 \longrightarrow\Lambda_{\ol{\mathbf b}}$
be a graph isomorphism, where
$\mathbf{a},\mathbf{b}\in\Ftwo^N$.
We apply Theorem~\ref{th:thfc} to the ordered fibre products $\Bipprod_{k=1}^{N}C_k^{(a_k)}$ and $\Bipprod_{k=1}^{N}C_k^{(b_k)}$, respectively.
There exist a permutation $\sigma\in S_N$, factor isomorphisms $f_k:\Und(C_k)\to\Und(C_{\sigma(k)})$, and a common sign $\delta\in\Ftwo$ such that for any
$z\in V(C_k)$ and $1\le k\le N$,  \begin{equation}\label{eq:eqth16pss}
 \mrc_{C_{\sigma(k)}}(f_k(z))+b_{\sigma(k)}
   =\mrc_{C_k}(z)+a_k+\delta,
\end{equation}
and, for any ${\bf u}\in V(\Lambda_{\ol{{\bf a}}})$ and $1\leq k\leq N$, we  have
\begin{equation}\label{eq:eqth16psco}
 [\varphi({\bf u})]_{\sigma(k)}=f_k(u_k).
\end{equation}

Define $w\in\widetilde{\mfw}$ on the entire refined direct product by the same $\sigma$ and $f_k$. By Eq.~\eqref{eq:eqth16psco}, $w({\bf u})$ and $\varphi({\bf u})$ have the same coordinate in each output position $\sigma(k)$. Since $\sigma$ permutes all coordinate positions, we have $w|_{\Lambda_{\ol{{\bf a}}}}=\varphi$.
Let $s_k$ denote the sign of $f_k$ as in Eq.~\eqref{eq:eqth16ps}. By Eqs.\eqref{eq:eqth16ps} and~\eqref{eq:eqth16pss}, cancellation of $\mrc_{C_k}(z)$ in $\Ftwo$ yields $b_{\sigma(k)}=a_k+s_k+\delta$ for $1\le k\le N$. Thus ${\bf b}=\sigma_*({\bf a}+{\bf s})+\delta\mathbf1_N$. Passing to the quotient, we obtain
\begin{equation*}
 \ol{{\bf b}}=\overline{\sigma_*({\bf a}+{\bf s})}=w\cdot\ol{{\bf a}}.
\end{equation*}
This proves that isomorphic components have indices in the same orbit. Conversely, if $\ol{{\bf b}}=w\cdot\ol{{\bf a}}$ for some $w\in\widetilde{\mfw}$, then by Eq.~\eqref{eq:eqth16pca}, we have $w(\Lambda_{\ol{{\bf a}}})=\Lambda_{\ol{{\bf b}}}$. Since $w$ is a graph automorphism, its restriction is an isomorphism between those components.

Finally, by part (\ref{it:ierfaa}), the maps $\Theta_{\ol{\boldsymbol{\epsilon}}}$ and $\Theta_{\ol{\boldsymbol{\eta}}}$ are isomorphisms onto $\Lambda_{\ol{d}(\ol{\boldsymbol{\epsilon}})}$ and $\Lambda_{\ol{d}(\ol{\boldsymbol{\eta}})}$, respectively. If $\psi:\Omega_{\ol{\boldsymbol{\epsilon}}}\to\Omega_{\ol{\boldsymbol{\eta}}}$ is an isomorphism, then
\[
 \Theta_{\ol{\boldsymbol{\eta}}}\circ\psi\circ\Theta_{\ol{\boldsymbol{\epsilon}}}^{-1}
 :\Lambda_{\ol{d}(\ol{\boldsymbol{\epsilon}})}\longrightarrow
  \Lambda_{\ol{d}(\ol{\boldsymbol{\eta}})}
\]
is an isomorphism, so  $\ol{d}(\ol{\boldsymbol{\epsilon}})$ and $\ol{d}(\ol{\boldsymbol{\eta}})$ lie in the same orbit. Conversely, if $\ol{d}(\ol{\boldsymbol{\epsilon}})$ and $\ol{d}(\ol{\boldsymbol{\eta}})$ lie in the same orbit, choose $w$ with $w\cdot \ol{d}(\ol{\boldsymbol{\epsilon}})=\ol{d}(\ol{\boldsymbol{\eta}})$. Then $\Theta_{\ol{\boldsymbol{\eta}}}^{-1}\circ   (w|_{\Lambda_{\ol{d}(\ol{\boldsymbol{\epsilon}})}})\circ\Theta_{\ol{\boldsymbol{\epsilon}}}$
 is an isomorphism from $\Omega_{\ol{\boldsymbol{\epsilon}}}$ to $\Omega_{\ol{\boldsymbol{\eta}}}$.
 This proves part~(\ref{it:ierfab}).

(\ref{it:ierfac}). Fix $\overline{\mathbf a}=\overline d(\ole)$.
By Eq.~\eqref{eq:eqth16pca}, restriction defines a homomorphism
\[
 \rse_{\overline{\mathbf a}}:
 \ostab_{\widetilde\mfw}(\overline{\mathbf a})
       \longrightarrow\aut(\Lambda_{\overline{\mathbf a}}),
 \]
given by $\rse_{\overline{\mathbf{a}}}(w)=w|_{V(\Lambda_{\overline{\mathbf a}})}$.
To prove surjectivity, take any $h\in\aut(\Lambda_{\overline{\mathbf a}})$. Apply the isomorphism construction in part~(\ref{it:ierfab}) with both source and target label equal to $\overline{\mathbf a}$ and with $\varphi=h$. It produces $w\in\widetilde{\mfw}$ whose restriction is exactly $h$, and whose label action satisfies $w\cdot\overline{\mathbf a}=\overline{\mathbf a}$. Thus $w$ lies in the stated stabilizer and
$\rse_{\overline{\mathbf a}}(w)=h$.

To prove injectivity, suppose $w=w_{\sigma,\mathbf f}$ fixes every vertex of this component. Each class of each $C_k$ has at least two vertices, by Lemma~\ref{le:leosdo}. Fix a global colour $t$. The vertices with that colour are the unrestricted Cartesian tuples $\prod_j C_j^{a_j+t}$. If $\sigma(k)\ne k$, choose two such tuples $\mathbf u,\mathbf u'$ differing only in position $\sigma(k)$.
They have $u_k=u_k'$ but $u_{\sigma(k)}\ne u_{\sigma(k)}'$. Since both are fixed,
\[
 u_{\sigma(k)}=[w({\bf u})]_{\sigma(k)}=f_k(u_k)
   =f_k(u'_k)=[w({\bf u}')]_{\sigma(k)}=u'_{\sigma(k)},
\]
a contradiction. We conclude that $\sigma(k)=k$ for every $k$, so $\sigma$ is the identity. If $N=1$, this conclusion is automatic, and the rest of the argument remains valid.

Now fix a position $k$ and a vertex $z\in V(C_k)$. Set $t=\mrc_{C_k}(z)+a_k$. For each $j\ne k$, since both colour classes of $C_j$ are nonempty, we can  choose a vertex $u_j$ of colour $a_j+t$ in $C_j$. Set $u_k=z$. Then ${\bf u}=(u_1,\ldots,u_N)\in V(\Lambda_{\ol{{\bf a}}})$. By $w({\bf u})={\bf u}$ and the identity $\sigma=1$, we have $f_k(z)=[w({\bf u})]_k=u_k=z$.
Since $z$ is arbitrary, $f_k$ is the identity on all of $V(C_k)$. This holds for every $k$, so every coordinate map is the identity. Consequently $w$ is the identity on all of $V(\widetilde{P})$. Thus $\ker\rse_{\ol{{\bf a}}}=1$, and $\rse_{\ol{{\bf a}}}$ is injective. Combining this with surjectivity, we have $\ostab_{\widetilde{\mfw}}(\ol{{\bf a}})\cong\aut(\Lambda_{\ol{{\bf a}}})$.

Finally, $\Theta_{\ol{\boldsymbol{\epsilon}}}$ is an isomorphism by part~(\ref{it:ierfaa}), so conjugation $g\mapsto\Theta_{\ol{\boldsymbol{\epsilon}}}^{-1}\circ g\circ\Theta_{\ol{\boldsymbol{\epsilon}}}$ is an isomorphism from $\aut(\Lambda_{\ol{{\bf a}}})$ to $\aut(\Omega_{\ol{\boldsymbol{\epsilon}}})$ for $g\in\aut(\Lambda_{\ol{{\bf a}}})$. Composing this isomorphism with $\rse_{\ol{{\bf a}}}$ gives the isomorphism asserted in part~(\ref{it:ierfac}).
This proves part~(\ref{it:ierfac}).

(\ref{it:ierfad}).
By part~(\ref{it:ierfaa}), we may identify $P$ with the
disjoint union of component graphs $\mathcal{P}:=\sqcup_{\ol{\mathbf{a}}\in \mathcal L_{\mathcal I}}\Lambda_{\ol{\mathbf{a}}}$, where $\mathcal L_{\mathcal I}\subseteq Q_N$ is a set of labels, and every
$\Lambda_{\ol{\mathbf{a}}}$ is an entire connected component. For $\mathcal{O}\in\mathscr{O}_{\mathcal I}$, define the subgraph $\mathcal{P}_{\mathcal{O}} :=\sqcup_{\ol{\mathbf{a}}\in\mathcal{O}\cap \mathcal L_{\mathcal I}}\Lambda_{\ol{\mathbf{a}}}$.
The nonempty intersections $\mathcal O\cap \mathcal L_{\mathcal I}$ partition $\mathcal L_{\mathcal I}$ since distinct orbits in $Q_N$ are disjoint.  Therefore $\mathcal P=\bigsqcup_{\mathcal O\in\mathscr O_{\mathcal I}}
 \mathcal P_{\mathcal O}$.

By part~(\ref{it:ierfab}), two selected connected components are isomorphic if and only if their labels lie in the same $\widetilde{\mfw}$-orbit on $Q_N$.  Consequently the graphs $\mathcal{P}_{\mathcal{O}}$
consist of exactly
$n_{\mathcal{O}}=|\mathcal{O}\cap \mathcal L_{\mathcal I}|$
mutually  isomorphic connected components,  and graphs from distinct $\mathcal{O}$ contain no isomorphic components.
Every automorphism of $\mathcal{P}$ must preserve each
$\mathcal{P}_{\mathcal{O}}$  setwise. Conversely, automorphisms of these disjoint subgraphs may be chosen independently, yielding
\begin{equation}\label{eq:th16-directsum-aut}
 \aut(\mathcal P)
 \cong\prod_{\mathcal O\in\mathscr O_{\mathcal I}}
       \aut(\mathcal P_{\mathcal O}).
\end{equation}

Fix $\mathcal{O}\in\mathscr{O}_{\mathcal{I}}$ and abbreviate $n=n_{\mathcal{O}}$ and
$\Lambda'=\Lambda_{\ol{{\bf a}}_{\mathcal{O}}}$.
Write $G=\aut(\Lambda')$. Enumerate the connected components of $\mathcal{P}_{\mathcal{O}}$ as $\Lambda_1,\ldots,\Lambda_n$. By part (\ref{it:ierfab}), each $\Lambda_j$ is isomorphic to
$\Lambda'$.  Fix an isomorphism
$\iota_j:\Lambda'\longrightarrow\Lambda_j$ for every $j$.
For $\mathbf g=(g_1,\ldots,g_n)\in G^n$, $\tau\in S_n$ and
$\mathbf u\in V(\Lambda_j)$, define
\begin{equation}\label{eq:eqth16pwm}
 \Psi_{\mathbf g,\tau}(\mathbf u)
   =\iota_{\tau(j)}\circ g_{\tau(j)}\circ
      \iota_j^{-1}(\mathbf u).
\end{equation}
This map is an automorphism that sends $\Lambda_j$ to
$\Lambda_{\tau(j)}$.  Conversely, every automorphism $\Psi$ of
$\mathcal P_{\mathcal O}$ permutes its connected components by a unique $\tau\in S_n$, and the restrictions  determine uniquely
$g_{\tau(j)}=\iota_{\tau(j)}^{-1}\circ
  (\Psi|_{\Lambda_j})\circ\iota_j$.
Thus $(\mathbf g,\tau)\mapsto\Psi_{\mathbf g,\tau}$ is bijective.

For ${\bf u}\in V(\Lambda_j)$, applying $\Psi_{{\bf h},\omega}$ and then $\Psi_{{\bf g},\tau}$ yields
\[
\begin{aligned}
 (\Psi_{{\bf g},\tau}\circ \Psi_{{\bf h},\omega})({\bf u})
 &=\iota_{\tau\omega(j)}\circ g_{\tau\omega(j)}
   \circ\iota_{\omega(j)}^{-1}
   \circ\iota_{\omega(j)}\circ h_{\omega(j)}\circ\iota_j^{-1}({\bf u})\\
 &=\iota_{\tau\omega(j)}\circ(g_{\tau\omega(j)}\circ h_{\omega(j)})
   \circ\iota_j^{-1}({\bf u}).
\end{aligned}
\]
Writing $\ell=\tau\omega(j)$, we have $\omega(j)=\tau^{-1}(\ell)$. Thus $\Psi_{{\bf g},\tau}\circ\Psi_{{\bf h},\omega}=\Psi_{{\bf g}(\tau\cdot{\bf h}),\tau\omega}$, where $(\tau\cdot{\bf h})_{\ell}=h_{\tau^{-1}(\ell)}$.
This is the multiplication in the semidirect product $G^n\rtimes S_n$ when $S_n$ acts on $G^n$ by $(\tau\cdot {\bf h})_\ell=h_{\tau^{-1}(\ell)}$. Consequently $\aut(\mathcal P_{\mathcal O})
 \cong G^n\rtimes S_n
 =G\wr S_n$.
The argument also applies when $n=1$ since $S_1$ is trivial. Together with  part (\ref{it:ierfab}), this identifies the factor in Eq.~\eqref{eq:eqgau} corresponding to this isomorphism class of components.

 Part~(\ref{it:ierfac}) further
identifies
$G\cong\ostab_{\widetilde{\mfw}}
       (\ol{\mathbf a}_{\mathcal O})$.
Substituting this into Eq.~\eqref{eq:th16-directsum-aut}
and transporting the resulting isomorphism through $\Theta$
gives Eq.~\eqref{eq:eqgau}, proving part~(\ref{it:ierfad}).

(\ref{it:ierfae}).
Use the type classes $J_\lambda$ and their sizes $t_\lambda$ from the statement. Fix an isomorphism
$\iota_k:\Und(C_k)\to T_{\lambda(k)}$, where $k\in J_{\lambda(k)}$. A factor symmetry has a unique coordinate permutation. Otherwise, suppose two representations used different input positions for one output position. Varying just one of those inputs would change that output in one representation but not the other. Each factor has at least two vertices, and its coordinate maps are injective, so this is impossible. With the permutation fixed, varying each input determines its map uniquely.

Consequently the map sending $w_{\sigma,\mathbf f}$ to $\sigma$ is a homomorphism onto $\Pi=\prod_{\lambda=1}^q S_{J_\lambda}$. Its kernel $\mathcal K$ is the coordinatewise group, which the maps $\iota_k$ identify with $\prod_{\lambda=1}^q\aut(T_\lambda)^{J_\lambda}$.
For $\sigma\in\Pi$, define a lift $\ell_\sigma$ by the factor maps $\iota_{\sigma(k)}^{-1}\iota_k$. Their composites satisfy
\[
 (\iota_{\tau\sigma(k)}^{-1}\iota_{\sigma(k)})
 (\iota_{\sigma(k)}^{-1}\iota_k)
       =\iota_{\tau\sigma(k)}^{-1}\iota_k,
\]
so $\ell_\tau\ell_\sigma=\ell_{\tau\sigma}$. They give a homomorphic section of the permutation map. In the transported coordinates, conjugation by $\ell_\sigma$ sends a tuple $\mathbf h$ in the kernel to the tuple with $j$-th entry $h_{\sigma^{-1}(j)}$. This is well defined since $\sigma$ preserves all $J_\lambda$. The semidirect product is therefore
\[
 \widetilde\mfw\cong
 \prod_{\lambda=1}^q
   \bigl(\aut(T_\lambda)^{J_\lambda}\rtimes S_{J_\lambda}\bigr)
 \cong\prod_{\lambda=1}^q
   \bigl(\aut(T_\lambda)\wr S_{t_\lambda}\bigr).
\]
This proves part (\ref{it:ierfae}). 
\end{proof}
\begin{remark}
The following observations clarify the choices in Theorem~\ref{th:thrfa}.
\begin{enumerate}
\item The description is independent, up to the indicated isomorphisms, of the chosen maps $\theta_i$. To see this, let $\theta'_i:B_i\to\Bipprod_{\ell\in I'_i}C'_\ell$ be another family of fibre-prime decompositions. Theorem~\ref{th:thci}, applied to $\theta'_i\theta_i^{-1}$, gives a bijection $\sigma_i:I_i\to I'_i$ and colour-preserving factor isomorphisms $f_k:C_k\to C'_{\sigma_i(k)}$. In particular the two refinements have the same total number $N$ of factors. Their union is a bijection
$\sigma$ of the two refined index sets and defines a coordinate isomorphism $F:\widetilde P\to\widetilde P'$ satisfying $F\circ\Theta=\Theta'$. Since every $f_k$ preserves colour, the induced bijection on labels is $\overline{\mathbf a}\mapsto \overline{\sigma_*\mathbf a}$. It maps $\mathcal L_{\mathcal I}$ onto
$\mathcal L_{\mathcal I'}$, because $\sigma(I_i)=I'_i$ for each original factor. Conjugation by $F$ sends the refined factor-symmetry group onto the other refined factor-symmetry group: a coordinate permutation and its
factor maps are carried to the corresponding conjugate coordinate data, and the inverse conjugation has the same property. The two label actions are therefore equivalent. Corresponding stabilizers are conjugate and the corresponding orbit intersections have equal
sizes. Hence changing factorization coordinates does not change the group or multiplicity data in Theorem~\ref{th:thrfa}.
\item
The label set $\mathcal L_{\mathcal I}$ need not be preserved by $\widetilde\mfw$. For example, take $B_1=C_6\bip C_6$ and $B_2=C_6$. Each $C_6$ is twin-free and fibre-prime by Lemma~\ref{le:lecp}(\ref{it:itcpa}). With $I_1=\{1,2\}$ and $I_2=\{3\}$, the label $\overline{(0,0,0)}$ is admissible. A colour-reversing automorphism in the first refined coordinate sends it to $\overline{(1,0,0)}$, which is not admissible.
Thus the multiplicities in Theorem~\ref{th:thrfa} are intersections of ambient orbits with $\mathcal L_{\mathcal I}$; they are not orbits of an
action on $\mathcal L_{\mathcal I}$.
\end{enumerate}
\end{remark}

We next derive consequences of Theorem~\ref{th:thrfa}.

\begin{corollary}\label{co:corfa}
Let $r\geq 1$, and let $B_1,\ldots,B_r$ be connected twin-free nonunit $f$-prime ordered bipartite graphs.
Put $P=\prod_{i=1}^r\Und(B_i)$ and $\mfw=\mfw(B_1,\ldots,B_r)$.
For $\ol{\boldsymbol{\epsilon}}\in Q_r$, put $\Omega_{\ol{\boldsymbol{\epsilon}}}=\Und\bigl(B_1^{(\epsilon_1)}\bip\cdots\bip B_r^{(\epsilon_r)}\bigr)$.
Then the following statements hold.
\begin{enumerate}[(i)]
\item\label{it:itrfa2} Two components $\Omega_{\ol{\boldsymbol{\epsilon}}}$ and $\Omega_{\ol{\boldsymbol{\eta}}}$ are isomorphic if and only if $\ol{\boldsymbol{\epsilon}}$ and $\ol{\boldsymbol{\eta}}$ lie in the same orbit of $\mfw$ on $Q_r$.
\item\label{it:itrfa3} For each $\ol{\boldsymbol{\epsilon}}\in Q_r$,
$\aut(\Omega_{\ol{\boldsymbol{\epsilon}}})\cong
 \ostab_{\mfw}(\ol{\boldsymbol{\epsilon}})$.
\item\label{it:itrfa4} Write $Q_r/\mfw$ for the set of orbits of $\mfw$ on $Q_r$. If $\mathcal{O}$ runs through these orbits and one representative $\ol{\boldsymbol{\epsilon}}_{\mathcal{O}}\in\mathcal{O}$ is chosen for each orbit, then
\begin{equation*}
 \aut(P)\cong\prod_{\mathcal{O}\in Q_r/\mfw}
 (\ostab_{\mfw}(\ol{\boldsymbol{\epsilon}}_{\mathcal{O}})
   \wr S_{|\mathcal{O}|}).
\end{equation*}
\item\label{it:itrfa5}
Let $F_1,\ldots,F_q$ be representatives of the pairwise nonisomorphic underlying graphs occurring among $\Und(B_1),\ldots,\Und(B_r)$. Let $J_\lambda=\{i:\Und(B_i)\cong F_\lambda\}$ and $t_\lambda=|J_\lambda|$.
Then
\begin{equation}\label{eq:eqfsw}
 \mfw \cong \prod_{\lambda=1}^q (\aut(F_\lambda)\wr S_{t_\lambda}).
\end{equation}
\end{enumerate}
\end{corollary}

\begin{proof}
Since each $B_i$ is connected, twin-free, nonunit, and $f$-prime, choose the factorization consisting only of $B_i$, with $\theta_i=\id_{B_i}$ in ~\eqref{eq:eqcr}. Then $N=r$, $I_i=\{i\}$ and $C_i=B_i$.
For $\boldsymbol\epsilon\in\Ftwo^r$, the $i$-th coordinate of $d(\boldsymbol\epsilon)$ is then $\epsilon_i$, so both $d$ and $\ol{d}$ are identity maps. Every vector is constant on each singleton block. Hence $\mathcal L_{\mathcal I}=Q_r$.
Since each $\theta_i$ is the identity and each block contains one coordinate, the map $\Theta$ is also the identity. Thus, $\Lambda_{\ol{\boldsymbol{\epsilon}}}=\Omega_{\ol{\boldsymbol{\epsilon}}}$ and $\Theta_{\ol{\boldsymbol{\epsilon}}}=\id_{\Omega_{\ol{\boldsymbol{\epsilon}}}}$.
Moreover, $\widetilde{\mfw}=\mfw$ and their actions on $Q_r$ agree. Theorem~\ref{th:thrfa} therefore applies.
Parts (\ref{it:ierfab}), (\ref{it:ierfac}), (\ref{it:ierfad}) and (\ref{it:ierfae}) of Theorem~\ref{th:thrfa} establish parts (\ref{it:itrfa2}), (\ref{it:itrfa3}), (\ref{it:itrfa4}) and (\ref{it:itrfa5}), respectively.
\end{proof}
The next corollary treats pairwise nonisomorphic factors and permits one unit factor. Its notation is local to the corollary and the following remark.

\begin{corollary}\label{th:thadp}
Let $r\geq1$, and let $B_1,\ldots,B_r$ be connected twin-free ordered bipartite graphs with pairwise nonisomorphic underlying graphs. Assume each $B_i$ is either $K_2$ or a nonunit $f$-prime graph. Put $\mathbf B=(B_1,\ldots,B_r)$ and $P=\prod_i\Und(B_i)$, and define
\[
\begin{aligned}
 \mathcal G_{\mathbf B}&=\{(g_i)_i:\ g_i\in\aut(\Und(B_i)),\
                \chi_{B_1}(g_1)=\cdots=\chi_{B_r}(g_r)\},\\
 L_{\mathbf B}&=\langle\overline{\mathbf e_i}:\autm(B_i)\ne\varnothing\rangle
      \leq Q_r,
\end{aligned}
\]
where $\mathbf e_i$ is the $i$-th standard basis vector and an empty
set of generators gives $L_{\mathbf B}=0$. Then the following hold.
\begin{enumerate}[(i)]
\item\label{it:itadp2}
For every $\ole\in Q_r$, $\aut(\Omega_{\ole})\cong\mathcal G_{\mathbf B}$.
\item\label{it:itadp3}
Two components $\Omega_{\ole}$ and $\Omega_{\overline{\boldsymbol\eta}}$
are isomorphic if and only if
$\overline{\boldsymbol\eta}-\ole\in L_{\mathbf B}$.
\item\label{it:itadp4}
Each coset of $L_{\mathbf B}$ indexes one isomorphism class of $|L_{\mathbf B}|$ components, and
\begin{equation}\label{eq:eqfap}
 \aut(P)\cong\prod_{\mathcal O\in Q_r/L_{\mathbf B}}
                  \bigl(\mathcal G_{\mathbf B}\wr S_{|L_{\mathbf B}|}\bigr).
\end{equation}
\end{enumerate}
\end{corollary}

\begin{proof}
We first assume that no factor is $K_2$. Since every factor is connected, twin-free, nonunit, and $f$-prime, Corollary~\ref{co:corfa} applies. If $w_{\sigma,\mathbf f}\in\mfw$, then $\Und(B_i)\cong\Und(B_{\sigma(i)})$ for every $i$. By pairwise nonisomorphism, $\sigma(i)=i$. Conversely, arbitrary automorphisms of the individual factors act coordinatewise. Thus, as a group of permutations of $V(P)$,
\[
 \mfw=\prod_{i=1}^r\aut(\Und(B_i)).
\]
For  $\mathbf g=(g_1,\ldots,g_r)$, define its colour sign vector by $\mathbf{s}(\mathbf{g})=(\chi_{B_i}(g_i))_{i=1}^r\in\Ftwo^r$. This map  is a homomorphism into the additive group $\Ftwo^r$. By Eq.~\eqref{eq:eqfsa}, $\mathbf{g}\cdot\overline{\boldsymbol\epsilon}=\overline{\boldsymbol\epsilon}+\overline{\mathbf s(\mathbf g)}$ .

(\ref{it:itadp2}). A tuple ${\bf g}$ fixes $\ol{\boldsymbol{\epsilon}}$ if and only if $\ol{{\bf s}({\bf g})}=0$. Since the kernel of $\Ftwo^r\to Q_r$ is $\langle\mathbf1_r\rangle$, $\ol{{\bf s}({\bf g})}=0$ means ${\bf s}({\bf g})=0$ or ${\bf s}({\bf g})=\mathbf{1}_r$. These alternatives are equivalent to $g_i\in \autp(B_i)$ for every $i$, or $g_i\in \autm(B_i)$ for every $i$, respectively. Thus $ \ostab_\mfw(\ol{\boldsymbol{\epsilon}})=(\prod_i\autp(B_i))\cup(\prod_i\autm(B_i))=\mathcal G_{\mathbf B}$.
Since $\ostab_\mfw(\ol{\boldsymbol{\epsilon}})$ is the kernel of ${\bf g}\mapsto\ol{{\bf s}({\bf g})}$, the indicated union is a subgroup. If some $\autm(B_i)$ is empty, the second alternative is impossible. By Corollary~\ref{co:corfa}(\ref{it:itrfa3}), restriction induces the asserted isomorphism.

If one factor, say $B_r$, is $K_2$ and $r\geq2$, deleting its coordinate identifies each component with the corresponding fibre product of the remaining factors. The preceding argument determines its automorphism group. Each such automorphism has a common coordinate sign and extends uniquely to the omitted coordinate by the automorphism of $K_2$ with that sign. This again identifies the component automorphism group with $\mathcal G_{\mathbf B}$.
If $r=1$ and $B_1=K_2$, the assertion follows from $\mathcal G_{\mathbf B}=\aut(K_2)$.

(\ref{it:itadp3}).
Suppose first that
\[
 \varphi:\Omega_{\ol{\boldsymbol{\epsilon}}}\longrightarrow\Omega_{\ol{\boldsymbol{\eta}}}
\]
is an ordinary isomorphism.

If no factor is $K_2$, then Corollary~\ref{cor:corsfc} and the pairwise nonisomorphism of the underlying factors imply that there are graph automorphisms $g_i\in\aut(\Und(B_i))$ and one sign $\delta\in\Ftwo$ such that for every $1\leq i\leq r$
\begin{equation}\label{eq:eqpcs}
 \eta_i=\epsilon_i+\chi_{B_i}(g_i)+\delta.
\end{equation}

If one factor, say $B_r$,  is $K_2$, first omit it from both fibre products and apply Corollary~\ref{cor:corsfc} to the remaining nonunit factors. By pairwise nonisomorphism, the factor permutation must be the identity.
 Thus, when $r\geq 2$, there are automorphisms $g_i\in\aut(\Und(B_i))$ for $1\leq i\leq r-1$ and a sign $\delta\in\Ftwo$ such that $\eta_i=\epsilon_i+\chi_{B_i}(g_i)+\delta$. In the unit case $r=1$, choose $\delta=0$. No nonunit coordinates occur in this case. Define $s_r=\eta_r+\epsilon_r+\delta\in\Ftwo$.
Since $K_2$ has a unique automorphism of each sign, choose $g_r\in\aut(K_2)$ with $\chi_{K_2}(g_r)=s_r$. Then Eq.~\eqref{eq:eqpcs} also holds for $i=r$.

Let ${\bf s}=(\chi_{B_1}(g_1),\ldots,\chi_{B_r}(g_r))$. By Eq.~\eqref{eq:eqpcs},
$\boldsymbol{\eta}-\boldsymbol{\epsilon}={\bf s}+\delta{\bf 1}_r$. Passing to the quotient by $\langle {\bf 1}_r\rangle$, we obtain $\ol{\boldsymbol{\eta}}-\ol{\boldsymbol{\epsilon}}=\ol{{\bf s}}$. The coordinate $s_i$ is $1$ only if $\autm(B_i)\neq\varnothing$, so $\ol{{\bf s}}\in L_{\mathbf B}$. This proves the necessity.

Conversely, assume $\ol{\boldsymbol{\eta}}-\ol{\boldsymbol{\epsilon}}\in L_{\mathbf B}$. By the definition of $L_{\mathbf B}$, choose a vector ${\bf s}=(s_1,\ldots,s_r)\in\Ftwo^r$ representing this difference such that $s_i=1$ only for indices $i$ with $\autm(B_i)\neq\varnothing$. For each $i$, choose $g_i\in \autp(B_i)$ if $s_i=0$ and $g_i\in \autm(B_i)$ otherwise.
The coordinatewise map ${\bf g}=(g_1,\ldots,g_r)$ is an automorphism of $P$. A vertex ${\bf x}=(x_1,\ldots,x_r)$ of $\Omega_{\ol{\boldsymbol{\epsilon}}}$ has coordinate colour vector $\boldsymbol{\epsilon}+t{\bf 1}_r$ for some $t\in\Ftwo$. The vertex ${\bf g}({\bf x})$ has coordinate colour vector $\boldsymbol{\epsilon}+{\bf s}+t{\bf 1}_r$, whose class in $Q_r$ is $\ol{\boldsymbol{\eta}}$. Hence the restriction of ${\bf g}$ is an ordinary isomorphism from $\Omega_{\ol{\boldsymbol{\epsilon}}}$ onto $\Omega_{\ol{\boldsymbol{\eta}}}$. This proves part~(\ref{it:itadp3}).

(\ref{it:itadp4}).
By part~(\ref{it:itadp3}), the isomorphism class containing $\Omega_{\ol{\boldsymbol{\epsilon}}}$ is indexed by the coset $\ol{\boldsymbol{\epsilon}}+L_{\mathbf B}$. Since every coset has $|L_{\mathbf B}|$ elements, this isomorphism class consists of $|L_{\mathbf B}|$ connected components. By part~(\ref{it:itadp2}), the automorphism group of each of these components is $\mathcal G_{\mathbf B}$. For a graph consisting of $m$ disjoint copies of  a connected graph $\Sigma$, fix an isomorphism from each copy to $\Sigma$. An automorphism may then act independently by an element of $\aut(\Sigma)$ on every copy and may also permute the $m$ copies. The chosen isomorphisms realize these permutations as automorphisms, and conjugation permutes the $m$ coordinate groups.
Hence its automorphism group is
\[
 \aut(\Sigma)^m\rtimes S_m=\aut(\Sigma)\wr S_m.
\]
Automorphisms cannot interchange nonisomorphic connected components. The full automorphism group is therefore the direct product of the groups associated with the different component isomorphism classes. Substituting $\Sigma=\Omega_{\ol{\boldsymbol{\epsilon}}}$ and $m=|L_{\mathbf B}|$, we obtain Eq.~\eqref{eq:eqfap}. This proves part~(\ref{it:itadp4}) and completes the proof.
\end{proof}

\begin{remark}\label{re:reas}
Under the hypotheses of Corollary~\ref{th:thadp}, if every factor has a
colour-reversing automorphism, then $L_{\mathbf B}=Q_r$. All $2^{r-1}$ components
are isomorphic and $\aut(P)\cong\mathcal G_{\mathbf B}\wr S_{2^{r-1}}$.
In general $\autm(B_i)$ is a coset, not a subgroup. Neither this formula
nor the definition of $\mathcal G_{\mathbf B}$ assumes a colour-reversing involution
or a splitting of the individual sign homomorphisms.
\end{remark}


\subsection{Perducts of fibre powers.}
We now apply Theorem~\ref{th:thrfa} to products of fibre powers of a single $f$-prime graph. The answer depends on whether this graph admits a colour-reversing automorphism. 
For an ordered bipartite graph $B$ and an integer $t\geq1$, write
$B^{\star t}=\Bipprod_{k=1}^tB$. Given positive integers
$m_1,\ldots,m_r$, with $r\geq1$, put
\[
 \mathbf m=(m_1,\ldots,m_r),\qquad
 N=\sum_{i=1}^r m_i,\qquad
 P_{\mathbf m}=\prod_{i=1}^r\Und(B^{\star m_i}).
\]
For each $1\leq j\leq \lfloor N/2\rfloor$, define
\begin{equation}\label{eq:eqblpo}
 \nu_j=
 \begin{cases}
 [z^j]\displaystyle\prod_{i=1}^r(1+z^{m_i}),&2j<N,\\[1mm]
 \tfrac12[z^j]\displaystyle\prod_{i=1}^r(1+z^{m_i}),&2j=N,
 \end{cases}
\end{equation}
where $[z^j]$ denotes the coefficient of $z^j$. In particular, $\nu_0=1$ since every $m_i$ is positive.

\begin{corollary}\label{co:corp}
Let $B$ be a connected twin-free nonunit $f$-prime ordered bipartite graph, and set $A_B=\aut(\Und(B))$. For $r\geq1$ and positive $m_1,\ldots,m_r$, use the notation above.
\begin{enumerate}[(i)]
\item\label{it:itrp2}
If $\autm(B)\ne\varnothing$, the $2^{r-1}$ connected components of
$P_{\mathbf m}$ are all isomorphic to $\Und(B^{\star N})$, and
\begin{equation}\label{eq:eqblr}
 \aut(P_{\mathbf m})\cong
 \left[\bigl(\autp(B)^N\cup\autm(B)^N\bigr)\rtimes S_N\right]
       \wr S_{2^{r-1}}.
\end{equation}
Here the union consists of tuples whose coordinate maps have a common
colour sign, and $S_N$ acts by permuting the coordinates.

\item\label{it:itrp3}
If $\autm(B)=\varnothing$, the indices $j$ with $\nu_j>0$ parametrize all component isomorphism classes. Class $j$ has $\nu_j$ components, each with automorphism group
\begin{equation}\label{eq:eqblmu}
 \mathcal H_j=\begin{cases}
 (A_B\wr S_j)\times(A_B\wr S_{N-j}),&2j<N,\\[1mm]
 (A_B\wr S_j)\wr S_2,&2j=N.
 \end{cases}
\end{equation}
Consequently,
\begin{equation}\label{eq:eqwbg}
 \aut(P_{\mathbf m})\cong
 \prod_{\substack{0\leq j\leq\lfloor N/2\rfloor\\\nu_j>0}}
          \bigl(\mathcal H_j\wr S_{\nu_j}\bigr).
\end{equation}
\end{enumerate}
\end{corollary}

\begin{proof}
Each $B^{\star m_i}$ is connected and twin-free by
Lemmas~\ref{le:lebp} and~\ref{le:lebpa}. Its two colour classes have sizes $|B^0|^{m_i}$ and $|B^1|^{m_i}$, both at least two by Lemma~\ref{le:leosdo}, so it is nonunit. Refine the $i$-th original factor into $m_i$ copies of $B$ on a block $I_i$. Theorem~\ref{th:thrfa} applies, with
$\widetilde\mfw\cong A_B^N\rtimes S_N$. The selected components correspond to the $2^{r-1}$ labels in
$\mathcal L_{\mathcal I}$, represented by vectors constant on each block.

First consider a single fibre power. A factor symmetry fixes the zero
label in $Q_t$ exactly when its sign vector is $\mathbf0_t$ or
$\mathbf1_t$. Its coordinate maps therefore lie in
$\autp(B)^t\cup\autm(B)^t$. This is a subgroup, being the inverse image
of $\langle\mathbf1_t\rangle$ under the coordinate sign homomorphism.
Coordinate permutations normalize it and give a complement in the
stabilizer. Theorem~\ref{th:thrfa}(iii) gives
\begin{equation}\label{co:corcg}
 \aut(\Und(B^{\star t}))\cong
       \bigl(\autp(B)^t\cup\autm(B)^t\bigr)\rtimes S_t.
\end{equation}
This includes $t=1$; the second set is empty if $\autm(B)=\varnothing$.

(\ref{it:itrp2}). Choose $\gamma\in\autm(B)$. For any $\mathbf a\in\Ftwo^N$, apply $\gamma$ where $a_k=1$ and the identity elsewhere. Equation~\eqref{eq:eqfsa} sends the zero label to $\overline{\mathbf a}$, so the ambient action on $Q_N$ is transitive. All selected components are therefore isomorphic to the zero component $\Und(B^{\star N})$. Its group is given by Eq.~\eqref{co:corcg}; Theorem~\ref{th:thrfa}(\ref{it:ierfad}) assembles the $2^{r-1}$ copies and proves Eq.~\eqref{eq:eqblr}.

(\ref{it:itrp3}). Every factor automorphism now preserves colours, so
$A_B^N$ fixes all labels and the action is induced by $S_N$.
For a binary vector $\mathbf a$, let $|\mathbf a|$ be its Hamming weight.
The value $j=\min\{|\mathbf a|,N-|\mathbf a|\}$ is unchanged by
coordinate permutations and by replacing a representative with its
complement. Conversely, labels with the same value have representatives
of weight $j$, and a bijection of their supports extends to a coordinate
permutation. Thus this value classifies the ambient orbits.
For an original component with label $\ole$, its value is
\[
 j=\min\left\{\sum_{i=1}^r m_i\epsilon_i,
                   N-\sum_{i=1}^r m_i\epsilon_i\right\}.
\]

To count the selected labels in an orbit, choose the blocks on which an
admissible representative equals one. Its weight is the sum of their
sizes, so $[z^j]\prod_i(1+z^{m_i})$ counts those representatives.
Distinct subsets of blocks count separately even when some sizes agree.
A quotient label has two complementary representatives. If $2j<N$,
exactly one has weight $j$; if $2j=N$, both have that weight.
Complementation has no fixed binary vector, so the middle coefficient
is even. Thus the orbit intersection has size $\nu_j$, and these
intersections partition the selected labels. In particular
$\sum_j\nu_j=2^{r-1}$, including the zero-labelled type.

For the stabilizer, choose a representative of weight $j$ with support
$T\subseteq\{1,\ldots,N\}$. A coordinate permutation fixes its quotient
label exactly when it maps $T$ to $T$ or to $T^c$.
If $2j<N$, the unequal sizes exclude interchange. Independent
coordinate maps and permutations on the two sets give
$(A_B\wr S_j)\times(A_B\wr S_{N-j})$.
If $2j=N$, a fixed bijection $T\to T^c$, extended by its inverse,
gives a coordinate involution. The subgroup preserving both sets is
$(A_B\wr S_j)^2$, and conjugation by this involution exchanges its two
factors. Every set-interchanging symmetry is the involution times a
set-preserving symmetry. Hence the stabilizer is $(A_B\wr S_j)\wr S_2$.
This interchange concerns coordinate positions, not a colour reversal
of the individual factor $B$.

Theorem~\ref{th:thrfa}(\ref{it:ierfac})--(\ref{it:ierfad}) now gives the groups in Eq.~\eqref{eq:eqblmu} and the product in Eq.~\eqref{eq:eqwbg}. The calculation uses ambient orbits intersected with $\mathcal L_{\mathcal I}$. 
When $r=1$, only the type $j=0$ is selected, so the formula also agrees
with Eq.~\eqref{co:corcg} for a single fibre power.
\end{proof}

The direct product power $\Und(B)^{\times r}$ is obtained by taking
$m_1=\cdots=m_r=1$. Then $N=r$ and $\nu_j=\binom rj$ for $2j<r$, while $\nu_{r/2}=\frac12\binom r{r/2}$ when $r$ is even. In the no-reversal case all these multiplicities are
positive, so there are exactly $\lfloor r/2\rfloor+1$ component types.

The next example illustrates how unequal block sizes affect the multiplicities.
\begin{example}
Take $B=P_5$, with ordered colour classes of sizes $3$ and $2$. It is connected, twin-free, and four-cycle-free, so Lemma~\ref{le:lecp} makes it $f$-prime. The unequal class sizes exclude colour reversal. Its only automorphisms are the identity and path reflection, giving $\aut(\Und(B))\cong\zz_2$. For $\mathbf m=(1,1,2)$, we have
\[
 (1+z)^2(1+z^2)=1+2z+2z^2+2z^3+z^4.
\]
Thus $\nu_0=1$, $\nu_1=2$, and $\nu_2=1$. The corresponding groups are $\mathcal H_0=\aut(\Und(B))\wr S_4$,
$\mathcal H_1=\aut(\Und(B))\times(\aut(\Und(B))\wr S_3)$, and $\mathcal H_2=(\aut(\Und(B))\wr S_2)\wr S_2$. Consequently,
\[
 \aut(P_{(1,1,2)})\cong
 \mathcal H_0\times(\mathcal H_1\wr S_2)\times\mathcal H_2.
\]
\end{example}

\subsection{Stability of graph pairs}
The next result determines the full group behind the stability
equivalence. The equivalence is also a consequence of the
factor-loopless theorem of Wang and Gao~\cite{WG2026};
the relation between the hypotheses is explained below.
\begin{corollary}\label{pr:prmp}
Let $X$ be a finite connected nonbipartite twin-free graph such that  $\Dcov(X)$ is $f$-prime with its natural ordered bipartition. Let $B$ be a connected twin-free nonunit $f$-prime ordered bipartite graph, and suppose $\Und(B)\not\cong\Dcov(X)$. Let $\chi_B$ be the colour-sign homomorphism. Then
\begin{equation}\label{eq:eqmpg}
 \aut(X\times\Und(B)) \cong \autp(\Dcov(X))\rtimes_{\varepsilon_B}\aut(\Und(B)),
\end{equation}
where
\begin{equation*}
 \varepsilon_B(g)(h)
 =\kappa_X^{\chi_B(g)}h\kappa_X^{-\chi_B(g)}
\end{equation*}
for $g\in\aut(\Und(B))$ and $h\in\autp(\Dcov(X))$.
Moreover, $(X,\Und(B))$ is stable if and only if $X$ is stable.
\end{corollary}

\begin{proof}
Fix the ordered bipartition $V(B)=B^0\sqcup B^1$. Recall that $\Dcov(X)=X\times K_2$ is the canonical double cover of $X$, with ordered bipartition $\Dcov(X)^0=V(X)\times\{0\}$ and $\Dcov(X)^1=V(X)\times\{1\}$.
Define
\begin{equation*}
\Phi:\Und(\Dcov(X)\bip B)\longrightarrow X\times \Und(B)
\end{equation*}
by $\Phi((x,i),b)=(x,b)$, for $b\in B^i$ and $(x,i)\in V(\Dcov(X))$.
We first verify that $\Phi$ is a graph isomorphism. For each $(x,b)\in V(X)\times V(B)$, the value $i=\mrc_B(b)$
is uniquely determined. Thus $((x,i),b)$ is the unique vertex of $\Dcov(X)\bip B$ mapped to $(x,b)$, so $\Phi$ is bijective.

Moreover, $((x,i),b)$ is adjacent to $((x',1-i),b')$
in the fibre product if and only if $\{x,x'\}\in E(X)$ and $\{b,b'\}\in E(B)$.
These are precisely the conditions for  $(x,b)$ and $(x',b')$ to be adjacent
in the ordinary direct product $X\times \Und(B)$. Thus $\Phi$ is a graph isomorphism.

Since $X$ is connected, twin-free, nonbipartite, $\Dcov(X)=X\times K_2$ is connected, bipartite and twin-free.
Thus both $\Dcov(X)$ and $B$ are connected, twin-free, ordered bipartite graphs. By hypothesis, both $\Dcov(X)$ and $B$ are $f$-prime and $\Dcov(X)\not\cong B$ as underlying graphs. By Corollary~\ref{th:thadp} and the isomorphism $\Phi$, we have
\begin{equation*}
\aut(X\times \Und(B))\cong \mathcal{G}(\Dcov(X),B).
\end{equation*}
Recall that $\mathcal{G}(\Dcov(X),B)=\{(a,g)\in
\aut(\Dcov(X))\times\aut(\Und(B)):\chi_{\Dcov(X)}(a)=\chi_B(g)\}$.
Consider the map
$
\phi:\mathcal{G}(\Dcov(X),B)\longrightarrow\aut(\Und(B))$ defined by $\phi(a,g)=g$ for $(a,g)\in \mathcal{G}(\Dcov(X),B)$.
This is a surjective group homomorphism. Indeed, for every $g\in\aut(\Und(B))$, the covering involution
$$
\kappa_X:(x,i)\mapsto(x,1-i)
$$
has colour sign $1$, so $\kappa_X^{\chi_B(g)}$ has the same colour sign as $g$. Hence $(\kappa_X^{\chi_B(g)},g)\in \mathcal{G}(\Dcov(X),B)$, and $\phi(\kappa_X^{\chi_B(g)},g)=g$.
The kernel of $\phi$ consists precisely of the pairs $(h,1)$ for which $\chi_{\Dcov(X)}(h)=0$. Thus $\ker(\phi)=\{(h,1):h\in\aut^{+}(\Dcov(X))\}$,
and hence $\ker(\phi)\cong\aut^{+}(\Dcov(X))$.
We therefore have an exact sequence
$$
1\longrightarrow
\aut^{+}(\Dcov(X))
\longrightarrow
\mathcal{G}(\Dcov(X),B)
\xrightarrow{\phi}
\aut(\Und(B))
\longrightarrow1.
$$
To show this exact sequence splits, define
\[\varsigma:\aut(\Und(B))\longrightarrow \mathcal{G}(\Dcov(X),B)\]
by $\varsigma(g)=(\kappa_X^{\chi_B(g)},g)$.
Since $\chi_B$ is a homomorphism and $\kappa_X^2=1$, we have
$\varsigma(gh)=\varsigma(g)\varsigma(h)$ for all $g,h\in\aut(\Und(B))$. Moreover, $
\phi(\varsigma(g))=g$. Thus $\varsigma$ is a section of $\phi$.
The complement acts on $\aut^{+}(\Dcov(X))$ by conjugation. For $h\in\aut^{+}(\Dcov(X))$ and $g\in\aut(\Und(B))$, we have
$$
\begin{aligned}
\varsigma(g)(h,1)\varsigma(g)^{-1}
&=
(\kappa_X^{\chi_B(g)},g)
(h,1)
(\kappa_X^{-\chi_B(g)},g^{-1})\\
&=
\left(
\kappa_X^{\chi_B(g)}
h
\kappa_X^{-\chi_B(g)},
1
\right).
\end{aligned}
$$
Thus the induced action
\begin{equation*}
\varepsilon_B(g)(h)
=\kappa_X^{\chi_B(g)}h\kappa_X^{-\chi_B(g)}.
\end{equation*}
Consequently, $
\aut(X\times \Und(B))\cong\aut^{+}(\Dcov(X))
\rtimes_{\varepsilon_B}\aut(\Und(B))$.

Now suppose that $X$ is stable. Then
$\aut(\Dcov(X))=\aut(X)\times\langle\kappa_X\rangle$,
where $\aut(X)$ is identified with its natural lift to
$\Dcov(X)$.  These lifted automorphisms preserve the two bipartition classes of $\Dcov(X)$, whereas $\kappa_X$ interchanges them. Therefore $\aut^{+}(\Dcov(X))=
\aut(X)$.

 Every lifted automorphism $\alpha$ satisfies $\kappa_X\alpha=\alpha\kappa_X$, so $\kappa_X$ centralizes
$\aut^{+}(\Dcov(X))$.
Consequently $\varepsilon_B(g)(h)=h$ for every $g\in\aut(\Und(B))$ and $h\in\aut^{+}(\Dcov(X))$.
Thus the semidirect product in  Eq.~\eqref{eq:eqmpg}
 is a direct product, and $
\aut(X\times \Und(B))\cong\aut(X)\times\aut(\Und(B))$.

It remains to identify this isomorphism with the natural coordinatewise action. Under $\Phi$, an element $\alpha\in\aut(X)$ acts by
\[
((x,i),b)\longmapsto((\alpha(x),i),b),
\]
which corresponds to $(x,b)\longmapsto(\alpha(x),b)$.
The element $\varsigma(g)=(\kappa_X^{\chi_B(g)},g)$ acts by $((x,i),b)\longmapsto((x,i+\chi_B(g)),g(b))$.
Since $\mrc_B(g(b))=\mrc_B(b)+\chi_B(g)$, the image lies in the fibre product, and $\Phi$ identifies this action with $(x,b)\longmapsto(x,g(b))$.
Thus each factor acts naturally on its corresponding coordinate.
Therefore, as permutation groups on $V(X)\times V(B)$,
$\aut(X\times \Und(B))=\aut(X)\times\aut(\Und(B))$.
Hence the graph pair $(X,\Und(B))$ is stable if $X$ is stable. Conversely, if $(X, \Und(B))$ is stable, by \cite[Theorem~1.3(a)]{WQX2025} we can obtain that $X$ is stable.
\end{proof}

\begin{remark}
The stability equivalence in Corollary~\ref{pr:prmp} is covered by~\cite[Theorem~1.3]{WG2026}. Here is the implication between the hypotheses. If $X\cong A\times C$ with both factors nonunit, Since $C$ is connected and nonbipartite, both $A$ and $C$ are connected and nonbipartite.
The map
\[
 ((a,i),(c,i))\longmapsto((a,c),i)
\]
is a colour-preserving isomorphism from $\Dcov(A)\bip\Dcov(C)$ to $\Dcov(A\times C)$.
It preserves the two edge conditions even when $A$ or $C$ has loops, and both double covers are nonunit. This contradicts the fact that $\Dcov$ is $f$-prime.

Also $X$ and $\Und(B)$ are directly coprime. Otherwise, Since $X$ is prime, $\Und(B)\cong X\times Y$ for a
connected bipartite graph $Y$. Order the classes of $Y$ so that this isomorphism preserves the classes of $B$.
Deleting the $K_2$-coordinate gives an isomorphism
$\Dcov(X)\bip Y\cong B$. The $f$-primality of $B$ forces
$Y\cong K_2$, contrary to $\Und(B)\not\cong\Dcov(X)$.
The pair is thus nontrivial and $X$ is factor-loopless in the terminology of~\cite{WG2026}. The preceding proof records a direct group calculation in this setting,
not an independent new stability theorem.
\end{remark}

\subsection{Automorphism groups of the canonical double covers of bipartite generalized Petersen graphs}\label{sec:gpg-application}
We now apply the preceding automorphism formulas to canonical double covers of bipartite generalized Petersen graphs, recovering the bipartite cases of \cite[Theorem~1.2(iii)]{QXZ2021}. Let $\kappa$ denote the unique nonidentity automorphism of $K_2$. We begin with the following lemma.

\begin{lemma}\label{co:comtk2}
Let $B$ be a connected twin-free nonunit fibre prime ordered bipartite graph. Then
\begin{equation*}
  \aut\bigl(\Und(B)\times K_2\bigr)
 \cong
 \aut(\Und(B))\wr S_2.
\end{equation*}
\end{lemma}

\begin{proof}
Equip $K_2$ with the ordered bipartition $K_2^0=\{0\}$ and $K_2^1=\{1\}$, and let $\kappa$ denote its unique nonidentity automorphism. The projection onto the first coordinate identifies $\Und(B\bip K_2)$ with $\Und(B)$. The ordinary direct product $\Und(B)\times K_2$ is the disjoint union of this component and $\Und(B\bip K_{2}^{\op})$, which is also isomorphic to $\Und(B)$.
We apply Corollary~\ref{th:thadp} to the two ordered bipartite factors $B_1=B$ and $B_2=K_2$.
Since $B$ is nonunit, $B\not\cong K_2$, so the two factors are nonisomorphic. Their sign sets are
$A_1^0=\aut^{+}(B)$, $A_1^1=\aut^{-}(B)$ and $A_2^0=\{1\}$, $A_2^1=\{\kappa\}$.

The component index group in Corollary~\ref{th:thadp} is $
Q_2=\Ftwo^2/\langle(1,1)\rangle$,
which has order $2$. Since $K_2$ admits the colour interchanging automorphism $\kappa$, the second coordinate contributes the generator $\ol{e}_2$ to the subgroup $H$. As $\ol{e}_2$ is the unique nonzero element of $Q_{2}$, we obtain $H=Q_2$.
Consequently, the two connected components of $\Und(B)\times K_2$ belong to the same $H$-orbit and hence are isomorphic.
It remains to determine the automorphism group of one component. The same sign condition in Corollary~\ref{th:thadp} gives
$\mathcal{G}=(\aut^{+}(B)\times\{1\})\cup
(\aut^{-}(B)\times\{\kappa\})$.

Define
$\psi:\aut(\Und(B))\longrightarrow \mathcal{G}$
by  $\psi(g)=\bigl(g,\kappa^{\chi_B(g)}\bigr)$, where
$\chi_B(g)\in\Ftwo$ denotes the sign of $g$.
We show that $\psi$ is a group isomorphism. Since $\chi_B$ is a homomorphism, $\chi_B(gh)=\chi_B(g)+\chi_B(h)$, and therefore
$$
\begin{aligned}
\psi(gh)
&=\bigl(gh,\kappa^{\chi_B(gh)}\bigr)\\
&=\bigl(gh,\kappa^{\chi_B(g)+\chi_B(h)}\bigr)\\
&=\bigl(g,\kappa^{\chi_B(g)}\bigr)
\bigl(h,\kappa^{\chi_B(h)}\bigr)\\
&=\psi(g)\psi(h).
\end{aligned}
$$
Thus $\psi$ is a homomorphism. It is injective since its first coordinate is $g$. It is surjective because every element of $\mathcal{G}$ has one of the two forms $(g,1)$, $g\in\aut^{+}(B)$, or $(g,\kappa)$, $g\in\aut^{-}(B)$.
Hence $\mathcal{G}\cong\aut(\Und(B))$.

Since $H=Q_2$ and $|Q_{2}|=2$, there is a single isomorphism class of connected components, consisting of the two components of $\Und(B)\times K_2$. The nontrivial element of $S_2$ interchanges these two components. Therefore, by Eq.~\eqref{eq:eqfap}, $
\aut(\Und(B)\times K_2)\cong \mathcal{G}\wr S_{2}=(\mathcal{G}\times \mathcal{G})\rtimes S_2$.
Using $\mathcal{G}\cong\aut(\Und(B))$, we obtain $
\aut(\Und(B)\times K_2)\cong\aut(\Und(B))\wr S_2$.
\end{proof}

Let $n$ and $k$ be integers with $n\geq3$ and $1\leq k<n/2$. The generalized Petersen graph $\gp(n,k)$ has vertex set $\{u_i,v_i:i\in\zz_n\}$
and edge set
\[
 \bigl\{\{u_i,u_{i+1}\},\{u_i,v_i\},\{v_i,v_{i+k}\}:i\in\zz_n\bigr\},
\]
where all subscripts are taken modulo $n$.

 Lemma~\ref{co:comtk2} now gives a formula for the automorphism group of the double cover of a bipartite generalized Petersen graph.
\begin{lemma}\label{th:thbgpgk2}
Let $n$ be even and $k$ odd, where $1\leq k<n/2$. Then
\begin{equation}\label{eq:eqgpgk2}
 \aut\bigl(\gp(n,k)\times K_2\bigr)
 \cong
 \aut(\gp(n,k))\wr S_2.
\end{equation}
\end{lemma}

\begin{proof}
Set $B=\gp(n,k)$.  By \cite[Proposition~4.3]{AL2009}, $B$ is bipartite.
Equip $B$ with either of its two ordered bipartitions. The graph $B$ is also connected, twin-free, and cubic.
By Lemma~\ref{le:lecp}, $B$ is a nonunit $f$-prime object. Lemma~\ref{co:comtk2} applies and gives Eq.~\eqref{eq:eqgpgk2}.
\end{proof}
Combining Lemma~\ref{th:thbgpgk2} with the classification of automorphism groups by Frucht, Graver and Watkins~\cite{FGW1971} gives the following corollary, which recovers the bipartite cases of \cite[Theorem~1.2(iii)]{QXZ2021}. We use the following notation.
\begin{align*}
 F(n)
 &=\langle\rho,\delta\mid
   \rho^n=\delta^2=1,\ \delta\rho\delta=\rho^{-1}\rangle,\\
 H(n,k)
 &=\langle\rho,\alpha\mid
   \rho^n=\alpha^4=1,\ \alpha\rho\alpha^{-1}=\rho^k\rangle,\\
 J(n,k)
 &=\langle\rho,\delta,\alpha\mid
   \rho^n=\delta^2=\alpha^2=1,\
   \delta\rho\delta=\rho^{-1},\
   \alpha\rho\alpha=\rho^k,\
   \alpha\delta\alpha=\delta\rangle.
\end{align*}
Here these group names are local to the generalized Petersen application. In particular $F(n)$ is the dihedral group of order $2n$.
The exceptional groups are specified by the presentations
in~\cite[Lemma~2.1]{QXZ2021}:
\[
\begin{aligned}
 E_8=\langle a,b,c\mid{}&
 a^8=b^2=c^3=1,\quad bab=a^{-1},\quad bcb=c^{-1},\\[-1mm]
 &cac=a^{-1},\quad ca^4=a^4c\rangle,\\
 E_{24}=\langle a,b,c\mid{}&
 (ca)^2=b^2=c^3=1,\quad bab=a^{-1},\\[-1mm]
 &bcb=c^{-1},\quad ca^4=a^4c\rangle.
\end{aligned}
\]
These presentations remove any choice of an unstated semidirect-product action. The presentation of $H(n,k)$ uses our right-to-left map convention: when $k^2\equiv-1\pmod n$, it is realized by the shift
$\rho(u_i)=u_{i+1}$, $\rho(v_i)=v_{i+1}$ and the automorphism
$\alpha(u_i)=v_{ki}$, $\alpha(v_i)=u_{ki}$.
All three edge types are preserved; $\alpha^2$ negates subscripts and
$\alpha\rho\alpha^{-1}=\rho^k$. The $4n$ maps $\rho^a\alpha^b$
($0\leq a<n$, $0\leq b<4$) are distinct, as their layer actions,subscript slopes, and then translations distinguish them. The relations reduce every word to this form. Identification with the full group in the generic case is the cited classification.

\begin{corollary}\label{co:coqxzbb}
Let $n$ be even and $k$ odd, where $1\leq k<n/2$, and write $A(n,k)=\aut\bigl(\Dcov(\gp(n,k))\bigr)
       =\aut\bigl(\gp(n,k)\times K_2\bigr)$.
If
\[
 (n,k)\notin
 \{(4,1),(8,3),(10,3),(12,5),(24,5)\},
\]
then
\[
 A(n,k)
 \cong
 \begin{cases}
 F(n)\wr S_2,
   & k^2\not\equiv\pm1\pmod n,\\[2mm]
 J(n,k)\wr S_2,
   & k^2\equiv1\pmod n,\\[2mm]
 H(n,k)\wr S_2,
   & k^2\equiv-1\pmod n.
 \end{cases}
\]
Moreover,
\begin{align*}
 A(4,1)&\cong(S_4\times\zz_2)\wr S_2,\\
 A(8,3)&\cong E_8\wr S_2,\\
 A(10,3)&\cong(S_5\times\zz_2)\wr S_2,\\
 A(12,5)&\cong(S_4\times S_3)\wr S_2,\\
 A(24,5)&\cong E_{24}\rtimes\zz_2\bigr)\wr S_2.
\end{align*}
\end{corollary}

\begin{proof}
Lemma~\ref{th:thbgpgk2} gives $A(n,k)\cong\aut(\gp(n,k))\wr S_2$. The classification in~\cite{FGW1971}, with the exceptional presentations recorded in~\cite[Lemma~2.1]{QXZ2021}, identifies the factor group as $F(n)$, $J(n,k)$, $H(n,k)$, or the listed exceptional group in precisely the stated cases. Substitution proves each formula. The classification of these factor groups is an external input; the double-cover step is
the preceding lemma.
\end{proof}

\vskip 0.2 cm
\noindent{\bf Acknowledgements}~~The author is supported by the National Natural Science Foundation of China (12401453). The author is grateful to Binzhou Xia for his feedback on the first draft of the paper.
\vskip 0.2 cm

\noindent{\bf Declaration on the use of AI}~~
ChatGPT was used for grammar checking, language polishing, and improving the clarity of the exposition.

\noindent{\bf Data Availability Statement}~~No data sets were generated or
analyzed in this study.

\vskip 0.2 cm
\noindent{\bf Conflict of Interest}~~
The author declares no conflicts of interest.

\end{document}